\documentclass[11pt,draftcls,onecolumn] {IEEEtran}

\renewcommand{\baselinestretch}{1.39}

\IEEEoverridecommandlockouts                              % This command is only
\usepackage{ifpdf}
\ifpdf
    \usepackage{graphicx}
    \usepackage{epstopdf}
    \def\Ebox#1#2{%
	\includegraphics[width=#1\hsize]{fig/#2}
    }
\else
    \usepackage{graphicx}
    \def\Ebox#1#2{%
	\includegraphics[width=#1\hsize]{fig/#2.eps}
    }
\fi

\usepackage{mathptmx} % assumes new font selection scheme installed
\usepackage{epsfig} % for postscript graphics files
\usepackage{mathtools}
\usepackage{times} % assumes new font selection scheme installed
\usepackage{amssymb}  % assumes amsmath package installed
\usepackage{amsmath} % assumes amsmath package installed
\usepackage{tabularx}
\usepackage{multirow}
\usepackage{algorithm}
\usepackage{algorithmic}
\usepackage[font=small,labelfont=bf]{caption}

\usepackage[dvipsnames,usenames]{color}
\usepackage[colorlinks,%
linkcolor=BrickRed,%
filecolor=BrickRed,%
citecolor=RoyalPurple,%
]{hyperref}
\usepackage{color}

\usepackage{tabularx}

\title{\LARGE \bf
Probabilistic Data Association-Feedback Particle Filter for Multiple Target Tracking Applications}

\author{Tao Yang and Prashant G. Mehta% <-this % stops a space
\thanks{T. Yang and P. G. Mehta are with the Coordinated Science Laboratory and
the Department of Mechanical Science and Engineering at the University of
Illinois at Urbana-Champaign (UIUC)
{\tt\small taoyang1@illinois.edu; mehtapg@illinois.edu}}% <-this % stops a space,
\thanks{Financial support from the AFOSR
grant FA9550-09-1-0190  and the NSF grant EECS-0925534 is gratefully acknowledged.}%
\thanks{The conference version of this paper appeared
in~\cite{taoyang_acc12,adamtilton_fusion12}.}%
}

\def\qed{\hspace*{\fill}~\IEEEQED\par\endtrivlist\unskip}%\mbox{\rule[0pt]{1.3ex}{1.3ex}}}

\def\Re{\mathbb{R}}
\def\R{\mathbb{R}}

\def\varble{\,\cdot\,}

\usepackage{eufrak}

\def\varble{\,\cdot\,}

\def\Inov{I} %So we can change symbol easily later

\def\Sec#1{Sec.~\ref{#1}}
\def\Fig#1{Fig.~\ref{#1}}

\def\notes#1{\marginpar{\tiny #1}\typeout{Notes!
Notes!
Notes!
}}
\renewcommand{\notes}[1]{\typeout{notes!}}

\newcommand{\field}[1]{\mathbb{#1}}

\def\Re{\field{R}}
\def\v{{\sf K}}

\def\Sec#1{Sec.~\ref{#1}}

\def\eqdef{\mathrel{:=}}

\def\clL{{\cal L}}

\def\Sec#1{Sec~\ref{#1}}

\def\Inov{I} %So we can change symbol easily later

\newcounter{rmnum}

\newenvironment{romannum}{\begin{list}{{\upshape (\roman{rmnum})}}{\usecounter{rmnum}
\setlength{\leftmargin}{6pt}
\setlength{\rightmargin}{4pt}
\setlength{\itemindent}{-1pt}
}}{\end{list}}

\newcounter{anum}

\newcommand{\ud}{\,\mathrm{d}}

\def\E{{\sf E}}
\def\Expect{{\sf E}}

\def\Prob{{\sf P}}

\def\Expect{{\sf E}}

\newcommand{\ith}{i^\text{th}}
\newcommand{\jth}{j^\text{th}}

\newcommand{\rth}{r^\text{th}}
\newcommand{\mth}{m^\text{th}}
\newcommand{\nth}{n^\text{th}}

\def\Mscr{\mathcal{M}}

\def\Nscr{\mathcal{M}}
\def\Pscr{\Pi}

\def\R{\mathbb{R}}

\def\Fig#1{Fig.~\ref{#1}}
\def\Sec#1{Sec.~\ref{#1}}

\def\UZ{\underline{\mathcal{Z}}}

\def\uX{\underline{X}}

\def\uZ{\underline{Z}}
\def\ua{\underline{\alpha}}
\def\uA{\underline{A}}

\def\ugamma{\underline{\gamma}}
\def\P{{\sf P}}
\def\v{{\sf K}}
\def\Inov{I}

\def\IEEEQEDclosed{\mbox{\rule[0pt]{1.3ex}{1.3ex}}}
\def\qed{\nobreak\hfill\IEEEQEDclosed}
\def\M{M}
\def\N{M}

\newcommand{\bmat}[1]{\begin{bmatrix}#1 \end{bmatrix}}

\newcommand{\vX}{X}

\newcommand{\bearing}{h}

\newcommand{\ProbAsso}{\beta}

\newcommand{\particleX}{X}

\def\varble{\,\cdot\,}

\def\eqdef{\mathrel{:=}}

\usepackage{ulem}

\newtheorem{theorem}{Theorem}[section]

\newtheorem{lemma}[theorem]{Lemma}

\newtheorem{example}{Example}
\newtheorem{remark}{Remark}

\begin{document}
\normalem
\maketitle

\vspace{-0.5in}
\begin{abstract}
This paper is concerned with the problem of tracking single or
multiple targets with multiple non-target specific observations (measurements).  For
such filtering problems with data association uncertainty, a novel feedback
control-based particle filter algorithm is introduced. The algorithm is
referred to as the {\em probabilistic data association-feedback particle
filter (PDA-FPF)}.  The proposed filter is shown to represent a
generalization -- to the nonlinear non-Gaussian case -- of the classical
Kalman filter-based probabilistic data association filter (PDAF). 
%to the general class of nonlinear non-Gaussian filtering problems. % Solution to the joint data association problem is referred to as the joint PDA-FPF (JPDA-FPF), which is a direct extension of the PDA-FPF.
% The PDA-FPF is based on the feedback particle filter (FPF)
% introduced in our earlier
% papers~\cite{taoyang_acc11,taoyang_cdc11,taoyang_cdc12}. The
One remarkable conclusion is that the proposed PDA-FPF
algorithm retains the innovation error-based feedback structure
of the classical PDAF algorithm, even in the nonlinear non-Gaussian case. 
% The modified filter prediction is a weighted sum of individual and population prediction with more confidence given to the population.
The theoretical
results are illustrated with the aid of numerical examples motivated by multiple target tracking applications.\end{abstract}

%%%%%%%%%%%%%%%%%%%%%%%%%%%%%%%%%%%%%%%%%%%%%%%%%%%%%%%%%%%%%%%%%%%%%%%%%%%%%%%
\section{Introduction} % 1 page
%%%%%%%%%%%%%%%%%%%%%%%%%%%%%%%%%%%%%%%%%%%%%%%%%%%%%%%%%%%%%%%%%%%%%%%%%%%%%%%
\label{sec:intro}

Filtering with data association uncertainty is important to a
number of applications, including, air and missile defense
systems, air traffic surveillance, weather surveillance, ground
mapping, geophysical surveys, remote sensing, autonomous
navigation and robotics~\cite{Bar-Shalom_IEEE_CSM,
prob_robotics_Thrun}.  For example, consider the
problem of multiple target tracking (MTT) with radar. The targets
can be multiple aircrafts in air defense, multiple ballistic
objects in missile defense, multiple weather cells in weather
surveillance, or multiple landmarks in autonomous navigation
and robotics.  In each of these applications, there exists data
association uncertainty in the sense that one can not assign, in an apriori manner,
individual observations (measurements) to individual targets.
% the objective is to
% track the individual target locations based on noisy observations
% obtained from the radar.

Given the large number of applications, algorithms for filtering problems with data
association uncertainty have been extensively studied in the past;
cf.,~\cite{Bar-Shalom_IEEE_CSM,Bar-Shalom_Proc_IEEE,Kyriakides_08,survey_data_assoc} and references therein. A typical
algorithm is comprised of two parts:
\begin{romannum}
\item A {\em filtering} algorithm
%(typically based on a Kalman filter)
  for tracking a single target, and
\item A {\em data association} algorithm for associating
  observations to targets.
\end{romannum}

Prior to mid-1990s, the primary tool for filtering was Kalman
filter or one of its extensions, e.g., extended Kalman filter.
The limitations of these tools in applications arise on account
of nonlinearities, both in the dynamic models (e.g.,
drag forces in ballistic targets) and in the sensor models (e.g., range or bearing).  The nonlinearities can lead
to a non-Gaussian multimodal conditional distribution.  For
such cases, Kalman and extended Kalman filters are known to
perform poorly; cf.,~\cite{Ristic_book_2004}. Since the advent
and wide-spread use of particle
filters~\cite{gorsalsmi93,DouFreGor01}, such filters are
becoming increasing relevant to single and multiple target tracking
applications; cf.,~\cite{Kyriakides_08, Ristic_book_2004,BlomBloem07,Gordon97,Avitzour95,BlomBloem11Fusion,BlomBloem2006,Hue00trackingmultiple_SP,Oh_09} and references
therein.
%TAO: Are there other survey type papers on PF for target tracking?
%
% However,
% particle filters are known to suffer from high variance and robustness
% issues, especially where
% dynamic models of target are uncertain and/or
% unstable~\cite{Daum_10,Bickel_08}.

The second part is the data association algorithm. The purpose
of the data association algorithm is to assign observations to
targets.  The complications arise due to multiple non-target
specific observations (due to multiple targets in the coverage
area), missing observations (probability of detection less than
one, e.g., due to target occlusion), false alarms (due to
clutter) and apriori unknown number of targets (that require
track initiation).

The earlier solutions to the data association problem considered assignments in a deterministic
manner: These include the simple but non-robust ``nearest
neighbor'' assignment algorithm and the multiple hypothesis
testing (MHT) algorithm, requiring exhaustive
enumeration~\cite{Reid79analgorithm}. However,
exhaustive enumeration leads to an NP-hard problem because
number of associations increases exponentially with time.

%\smallskip

The complexity issue led to development of probabilistic
approaches: These include the probabilistic MHT or its simpler
``single-scan'' version, the probabilistic data
association (PDA) filter for tracking a single target in clutter, 
or its extension, the joint PDA (JPDA) filter for tracking multiple 
targets~\cite{Bar-Shalom_Proc_IEEE, Bar-Shalom_book_88}.
These algorithms require 
computation (or approximation) of the {\em
observation-to-target association probability}.  Certain
modeling assumptions are necessary to compute these in a
tractable manner. 

For certain MTT applications, the JPDA filter algorithm was found to
coalesce neighboring tracks for closely spaced
targets~\cite{Fitzgerald85,BlomBloem2000}. In order to overcome the
track coalescence issue, several extensions of the basic JPDA
algorithm have been developed, including the JPDA*
filter~\cite{BlomBloem2000,BlomBloem95CDC} and the set JPDA (SJPDA)
filter~\cite{SvenssonWillett11}. In MTT applications involving
multiple maneuvering targets in the presence of clutter, the JPDA
algorithm is combined with the interacting multiple model (IMM)
approach,  which leads to the IMM-JPDA class of
filters~\cite{ChenTugnait01,BlomBloem02Fusion}. Another extension of
the JPDA filter for tracking distinctive targets appears in~\cite{MTT_Tomlin_acc04}.  Even though 
PDA algorithms have reduced computational complexity, one limitation
is that they have been developed primarily in linear settings and rely on
Gaussian approximation of the posterior distribution;
cf.,~\cite{Bar-Shalom_IEEE_CSM}. 

%\smallskip

The rapid development of particle filtering has naturally led to
an investigation of data association algorithms based on
importance sampling techniques.  This remains an active area of
research; cf.,~\cite{survey_data_assoc,Herman_2011} and references therein. Two classes of particle filters for MTT applications with data association
uncertainty appear in the literature:
\begin{romannum}
\item The first approach involves application of the standard particle
  filtering approach to the joint (multiple target) state. For the
  joint state space, sequential importance resampling (SIR) particle
  filter is developed in~\cite{Gordon97,Avitzour95}. Unlike the
  classical (linear) PDA algorithm, the likelihood of observation needs to be calculated
  for each particle. This calculation represents the main
  computational burden.  Several extensions have also been considered:
  For the problem of track coalescence avoidance, a decomposed particle filter is introduced in~\cite{BlomBloem11Fusion}.  For multiple maneuvering target tracking in clutter, the joint IMM-PDA particle filter appears in~\cite{BlomBloem2006,BlomBloem03Fusion,BlomBloem02CDC}, where the IMM step and the PDA step are performed jointly for all targets. 

\item The second approach involves evaluation of association
  probability by using a Markov chain Monte Carlo (MCMC) algorithm,
  e.g., by randomly sampling from the subset 
  %  of the probability
  % simplex
  where the posterior is concentrated.  This avoids the
  computationally 
  intensive enumeration of all possible associations.  One early
  contribution is the multi-target particle filter (MTPF) introduced
  in~\cite{Hue00trackingmultiple_SP}.  In MTPF, samples of the association
  random variable are obtained iteratively from their joint posterior
  via the Gibbs sampler. In~\cite{Oh_09}, the Markov chain Monte Carlo
  data association (MCMCDA) filter is introduced where the association
  probability is approximated by using the Metropolis-Hastings
  algorithm.  
\end{romannum}
%Another direction of particle filtering methods for data association problem is that: instead of enumerating over all possible association events, a Markov chain Monte Carlo (MCMC) sampling techniques is applied to obtain the association probabilities. 
%By using MCMC techniques, these approaches can avoid the combinatorial drawbacks exists in the PDA-based data association algorithms. 

%A particle filter algorithm for MTT applications with
%kinematic constraints appears in~\cite{Kyriakides_08}. 

In applications, the sampling-based approaches may suffer from some the drawbacks of
particle filtering, e.g., particle impoverishment and degeneracy, slow convergence rate
(particularly for MCMC), robustness issues including sample variance
and numerical instabilities~\cite{DaumHuang10,Bickel_08,doucet2000sequential}. Moreover, the feedback structure of the Kalman filter-based PDA filter is no longer preserved.

\smallskip

In this paper, we introduce a novel feedback control-based
particle filter algorithm
for solution of the joint filtering-data association
problem. The proposed algorithm
%, joint probability data association-feedback particle filter
%(JPDA-FPF),
is based on the feedback particle filter (FPF) concept introduced by
us in earlier
papers~\cite{taoyang_TAC12,taoyang_acc11,taoyang_cdc11,taoyang_cdc12}. 
FPF is a controlled system where
the state of each particle evolves according to,
\[
\frac{\ud }{\ud t} [\text{Particle}] = [\text{Dynamics}] + [\text{Control}],
\]
where
\[
[\text{Control}] = [\text{Gain}]\cdot [\text{Innovation error}],
\]
and
\begin{align*}
[\text{Innovation error}]  = [\text{Observation}]-  \left( \frac{1}{2}\,
  [\text{Part. predict.}] + \frac{1}{2} \,[\text{Pop. predict.}] \right).
\end{align*}
The terms ``[Part. predict.]'' and ``[Pop. predict.]'' refer to
the prediction -- regarding the next ``[Observation]'' (value) -- as made by the particle and
by the population, respectively.  This terminology is made precise in the following
Section, with details
in~\cite{taoyang_TAC12, taoyang_acc11,taoyang_cdc11,taoyang_cdc12}.  In these papers, it is also
shown that, in the infinite
particle limit, FPF provides an {\em exact} solution of the nonlinear filtering
task. 
  
\begin{figure*}
    \centering
%    \hspace{-0.5in}
 \includegraphics[scale=0.31]{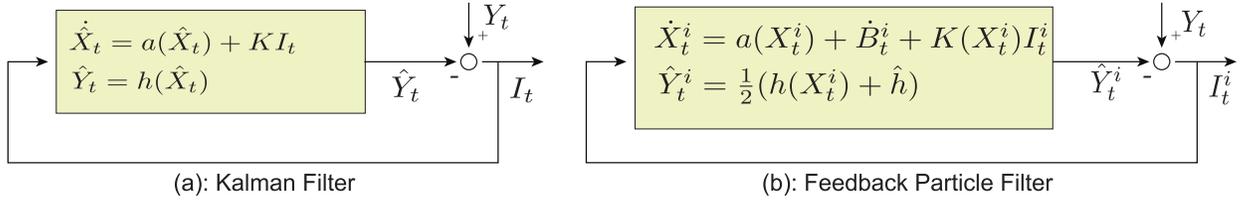}
   \vspace{-0.1in}
\caption{Innovation error-based feedback structure for (a) Kalman
  filter and (b) nonlinear feedback particle filter (see Remark~\ref{rem:rem1}).
}\vspace{-.4cm}
    \label{fig:fig_FPF_KF}
\end{figure*}

In contrast to a conventional particle filter, FPF does not
require resampling and thus does not suffer particle degeneracy.  FPF has a feedback structure similar to Kalman
filter (See Fig.~\ref{fig:fig_FPF_KF}).  Feedback is important on account of 
robustness:  Numerical results in~\cite{taoyang_TAC12,taoyang_cdc11} show that feedback can help
reduce the high variance that is sometimes observed in the
conventional particle filter.

\medskip

% The filter has a
% feedback structure similar to the Kalman filter: At each time
% $t$, the control is obtained by using a proportional gain
% feedback with respect to a certain modified form of the
% innovation error:
% \begin{equation}
% \text{Innov. error} = \text{Obs.} - \left[\frac{1}{2}\;\text{Ind. Pred.} + \frac{1}{2}\;\text{Pop. Pred.}\right],\nonumber
% \end{equation}
% where the coefficient $\frac{1}{2}$ indicates that the filter prediction gives equal weights on individuals as well as population, while in the classical nonlinear filtering equation, the prediction term is entirely due to the population. Then the particles are updated by using the proportional control term:
% \[
% \frac{\ud }{\ud t} [\text{Particle}] = [\text{Dynamics}] + \underbrace{[\text{Gain}] \cdot [\text{Innovation error}]}_{\text{Control input}}.
% \]
% The filter design amounts to design of the proportional gain -- the solution is given by the Kalman gain
% in the linear Gaussian case and obtained by solving a certain boundary value problem in the nonlinear case. The FPF has several advantages: since no resampling required in FPF, sampling related issue with conventional particle filter (i.e., particle degeneracy) are avoided; and the feedback structure can help improve robustness and reduce high variance.

In the present paper, we extend the basic feedback particle filter to
problems with data association uncertainty.  We refer to the
resulting algorithm as the probabilistic data
association-feedback particle filter (PDA-FPF).  As the name
suggests, the proposed algorithm represents a generalization of
the Kalman filter-based PDA filter now to the general nonlinear
non-Gaussian problems.  Just as with the classical PDA filter, an extension to multiple
target tracking case is also easily obtained.  The resulting algorithm
is referred to as JPDA-FPF.

One remarkable conclusion of our paper is that the PDA-FPF
retains the innovation error-based feedback structure even for
the nonlinear non-Gaussian problems.  The data association
uncertainty serves to modify the control:  Gain is decreased by a factor proportional to the
association probability, and the innovation error is based on a
modified expression for 
prediction whereby the particle gives additional weight to the population:
\begin{align*}
[\text{Gain}] & = \beta \cdot [\text{Nominal Gain}], \\
[\text{Innovation error}] &=[\text{Observation}] - \left(
  \frac{\beta}{2}\, [\text{Part. predict.}] +
  (1- \frac{\beta}{2})\, [\text{Pop. predict.}] \right),
\end{align*}    
where $\beta$ denotes the association probability (i.e., the
probability that the observation originates from the target). 
% Because of the clutter and unknown observation-to-target association, the modified innovation error turns out to have the following simple form:
% \begin{equation}
% \text{Innov. error} = \text{Obs.} - \left[\frac{\text{assoc. prob.}}{2}\; \text{Ind. Pred.} + \left(1 -\frac{\text{assoc. prob.}}{2}\right) \text{Pop. Pred.}\right],\nonumber
% \end{equation}
% which effectively means that the filter prediction puts more confidence on the population rather than individuals when the association is uncertain. This modified structure is exploited more in~\Sec{sec:pda_fpf_main}. 
The innovation error-based
feedback structure is expected to be useful because of the coupled
nature of the filtering and the data association problems.

The theoretical results are illustrated with the aid of three numerical
examples: i) tracking of a
single target in the presence of clutter, ii) tracking of two targets
in the presence of data association uncertainty and
iii) a multiple target tracking problem involving track coalescence.
Comparisons with the sequential importance resampling particle filter
(SIR-PF) are also provided.  
% For the track coalescence example, a model problem
% scenario from~\cite{BlomBloem2000} is used to illustrate possible track coalescence in the presence of data
% association uncertainty. A Monte Carlo simulation of 100 runs is conducted to compare the performance of JPDA-FPF to the SIR joint particle filter. %Finally a three-target tracking problem is examined.

The outline of the remainder of this paper is as follows:  We
begin with a brief review of the feedback particle filter in
\Sec{sec:prelim}.  The PDA-FPF algorithm is described
for the problem of tracking single target in the presence of clutter, in
\Sec{sec:PDAFPF}.  The JPDA-FPF for multiple target tracking case follows as a straightforward
extension, and is discussed in \Sec{sec:JPDAFPF}.  Numerical
examples appear in \Sec{sec:numerics}.

\section{Preliminaries: Feedback Particle Filter}
\label{sec:prelim}

%%Describe it in the multidimensional setting%%%%%%%%%%

In this section we briefly summarize the feedback particle
filter, introduced in our earlier
papers~\cite{taoyang_TAC12, taoyang_acc11, taoyang_cdc11, taoyang_cdc12}. In
these papers, we considered the following nonlinear filtering problem:
%\vspace{-15pt}
\begin{subequations}
\begin{align}
\ud X_t &= a(X_t)\ud t + \ud B_t,
\label{eqn:Signal_Process}
\\
\ud Z_t &= h(X_t)\ud t + \ud W_t,
\label{eqn:Obs_Process}
\end{align}
\end{subequations}
where $X_t\in\Re^d$ is the state at time $t$, $Z_t \in\Re^s$ is the
observation process, $a(\varble)$, $h(\varble)$ are $C^1$
functions, and $\{B_t\}$, $\{W_t\}$ are mutually independent
Wiener processes of appropriate dimensions.  The covariance matrix of
the observation noise $\{W_t\}$ is assumed to be positive definite. By scaling, we may assume, without loss of generality, that the covariance matrices associated with $\{B_t\}$, $\{W_t\}$ are identity matrices. The function $h$ is a column vector whose $\jth$ coordinate is denoted as $h_j$ (i.e., $h = (h_1,h_2,\hdots,h_s))$.

The objective of the filtering problem is to estimate the
posterior distribution, denoted as $p^*$, of $X_t$ given the history
$\mathcal{Z}_t := \sigma(Z_{\tau}:  \tau \le t)$. The evolution of
$p^*(x,t)$ is described by the Kushner-Stratonovich (K-S)
equation (see~\cite{jazwinski70}):
\begin{equation}
\ud p^\ast = \clL^\dagger p^\ast \ud t + ( h-\hat{h} )^T (\ud Z_t - \hat{h} \ud t)p^\ast, \nonumber
%\label{eqn:prelim_kushner}
\end{equation}
where
\begin{equation}
 \clL^\dagger p^\ast = - \nabla \cdot(p^\ast a) + \frac{1}{2}\Delta p^\ast, \label{eqn:dagger}
\end{equation}
and $ \hat{h} = \int_{\Re^d} h(x) p^*(x,t) \ud x$. Here, $\Delta$ denotes the Laplacian in $\Re^d$. If $a(\varble)$,
$h(\varble)$ are linear functions, the solution is given by the
finite-dimensional Kalman filter.

The feedback particle filter is a controlled system comprising
of $N$ particles. The dynamics of the $i^{\text{th}}$ particle has the following Stratonovich form:
\begin{align}
\ud X_t^i = a(X_t^i) \ud t + \ud B_t^i + \v(X_t^i,t) \circ \ud I_t^i,\label{eqn:prelim_fpf}
\end{align}
where $\{B_t^i\}$ are mutually independent standard Wiener processes, $\Inov_t^i$ represents a modified form of the \textit{innovation process} that
appears in the nonlinear filter,
\begin{equation}
\ud \Inov^i_t \eqdef \ud Z_t - \frac{1}{2}\left(h(X^i_t) + \hat{h}\right)\ud t,
\label{e:in}
\end{equation}
where   $\hat{h} := {\sf E} [h(X^i_t)|\mathcal{Z}_t] = \int_{\Re^d}
h(x) p(x,t) \ud x$ and $p(x,t)$ denotes the conditional
distribution of $X_t^i$ given $\mathcal{Z}_t$. In a numerical
implementation, we approximate $\hat{h} \approx \frac{1}{N}
\sum_{i=1}^N h(X^i_t) $.

The gain function $\v(x,t)$ is obtained as a solution to an Euler-Lagrange boundary value problem (E-L BVP) based on $p$: For $j = 1, \hdots, s$, the function $\phi_j: \Re^d \rightarrow \Re$ is a solution to the second order differential equation:
\begin{equation}
\label{eqn:EL_phi_prelim}
\begin{aligned}
\nabla \cdot (p(x,t) \nabla \phi_j(x,t) ) & = - (h_j(x)-\hat{h}_j) p(x,t),\\
\int_{\Re^d} \phi_j(x,t) p(x,t) \ud x & = 0,
\end{aligned}
\end{equation}
In terms of these solutions, the gain function is given by
\begin{equation}
[\v]_{lj} = 
\frac{\partial \phi_j} {\partial x_l}, \;\;\text{for} \;\;l \in \{1,\hdots,d\},\;  j\in\{1,\hdots,s\}.
\label{eqn:gradient_gain_fn_prelim}
\end{equation}
Note that the gain function needs to be obtained for each value of time $t$.

The evolution of $p(x,t)$ is given by a forward Kolmogorov operator (See Proposition 1 in~\cite{taoyang_cdc12}). 

In~\cite{taoyang_cdc12}, it is shown that the FPF~\eqref{eqn:prelim_fpf}-\eqref{eqn:gradient_gain_fn_prelim} is consistent: That is, given $p^\ast (x,0) = p(x,0)$ and the gain function $\v(x,t)$ is obtained according
to~\eqref{eqn:EL_phi_prelim}-\eqref{eqn:gradient_gain_fn_prelim}, then $p(x,t) = p^\ast (x,t)$ for all $t \geq 0$. This means that the empirical distribution of particles approximates the true posterior $p^\ast$ as the number of particles $N\rightarrow \infty$.

%
%\begin{figure*}
%    \centering
%%    \hspace{-0.5in}
% \includegraphics[scale=0.31]{figures/FPF_KF_Discrete}
%   \vspace{-0.1in}
%\caption{Innovation error-based feedback structure for (a) Kalman
%  filter and (b) nonlinear feedback particle filter (see Remark~\ref{rem:rem1}).
%}\vspace{-.4cm}
%    \label{fig:fig_FPF_KF}
%\end{figure*}

\begin{remark}
\label{rem:rem1}
Given that the Stratonovich form provides a mathematical
interpretation of the (formal) ODE model \cite[Section 3.3 of the SDE
text by {\O}ksendal]{Oksendal_book}, we also obtain the (formal) ODE
model of the filtering problem and the filter. Denoting $Y_t \doteq \frac{\ud Z_t}{\ud t}$ and white noise
processes $\dot{B}_t \doteq \frac{\ud B_t}{\ud t} $, $\dot{W}_t \doteq \frac{\ud W_t}{\ud t} $, the equivalent ODE model of~\eqref{eqn:Signal_Process} and~\eqref{eqn:Obs_Process} is:
\begin{align*}
\frac{\ud X_t}{\ud t} &= a(X_t) + \dot{B}_t,
\\
Y_t &= h(X_t) + \dot{W}_t,
\end{align*}
where $\dot{B}_t$ and $\dot{W}_t$ are mutually independent white noise processes.
The corresponding ODE model of the
filter is given by,
\begin{equation}
\frac{\ud X^i_t}{\ud t} = a(X^i_t) + \dot{B}^i_t +  \v(X^i,t) \cdot \left( Y_t - \frac{1}{2}    (h(X^i_t) + \hat{h})\right), \nonumber
%\label{eqn:FPF_ODE}
\end{equation}
where $\dot{B}^i_t \doteq \frac{\ud B^i_t}{\ud t}$ denotes the standard white noise process.
The feedback particle filter thus provides a generalization of the
Kalman filter to nonlinear systems, where the innovation error-based
feedback structure of the control is preserved (see
Fig.~\ref{fig:fig_FPF_KF}).  For the linear case, it is shown in
\cite{taoyang_cdc11,taoyang_cdc12} that the gain function
is the Kalman gain.  For the nonlinear case, the Kalman gain is
replaced by a nonlinear function of the state, obtained as the
solution
of~\eqref{eqn:EL_phi_prelim}-\eqref{eqn:gradient_gain_fn_prelim}. Various
approximation methods to obtain the nonlinear gain function have been
proposed in our earlier
papers~\cite{taoyang_acc11,taoyang_cdc11,taoyang_cdc12}. The Galerkin
and constant gain approximations are briefly reviewed in the following
subsection~\ref{apdx:const_approx}. 

Considering the wide use of the ODE formalism in practice and its (formal) equivalency to the SDE model~\eqref{eqn:prelim_fpf}, we make the following convention throughout the rest of paper: The filter is first derived using the rigorous SDE formalism. Following that, an equivalent ODE model for the filter is included as a remark. 
\qed
\end{remark}

\subsection{Galerkin and constant gain approximations}
\label{apdx:const_approx}

% The main difficulty in implementing the FPF algorithm is to obtain the
% solution of the BVP~\eqref{eqn:EL_phi_prelim} at each time step.  In
% this section,  

For a fixed time $t$ and $j \in \{1,\hdots, s\}$, a vector-valued function $\nabla \phi_j(x,t)$ is said to be a weak solution of the BVP~\eqref{eqn:EL_phi_prelim} if
\begin{equation}
\E\left[\nabla \phi_j \cdot \nabla \psi \right] = \E[(h_j-\hat{h}_j)\psi]\label{eqn:weak_form}
\end{equation}
holds for all $\psi \in H^1(\Re^d;p)$ where $\E[\cdot] := \int_{\Re^d} \cdot p(x,t)\ud x$ and $H^1$ is a certain Sobolev space (see~\cite{taoyang_cdc12}). 

\begin{figure}
\centering
\includegraphics[scale = 0.5]{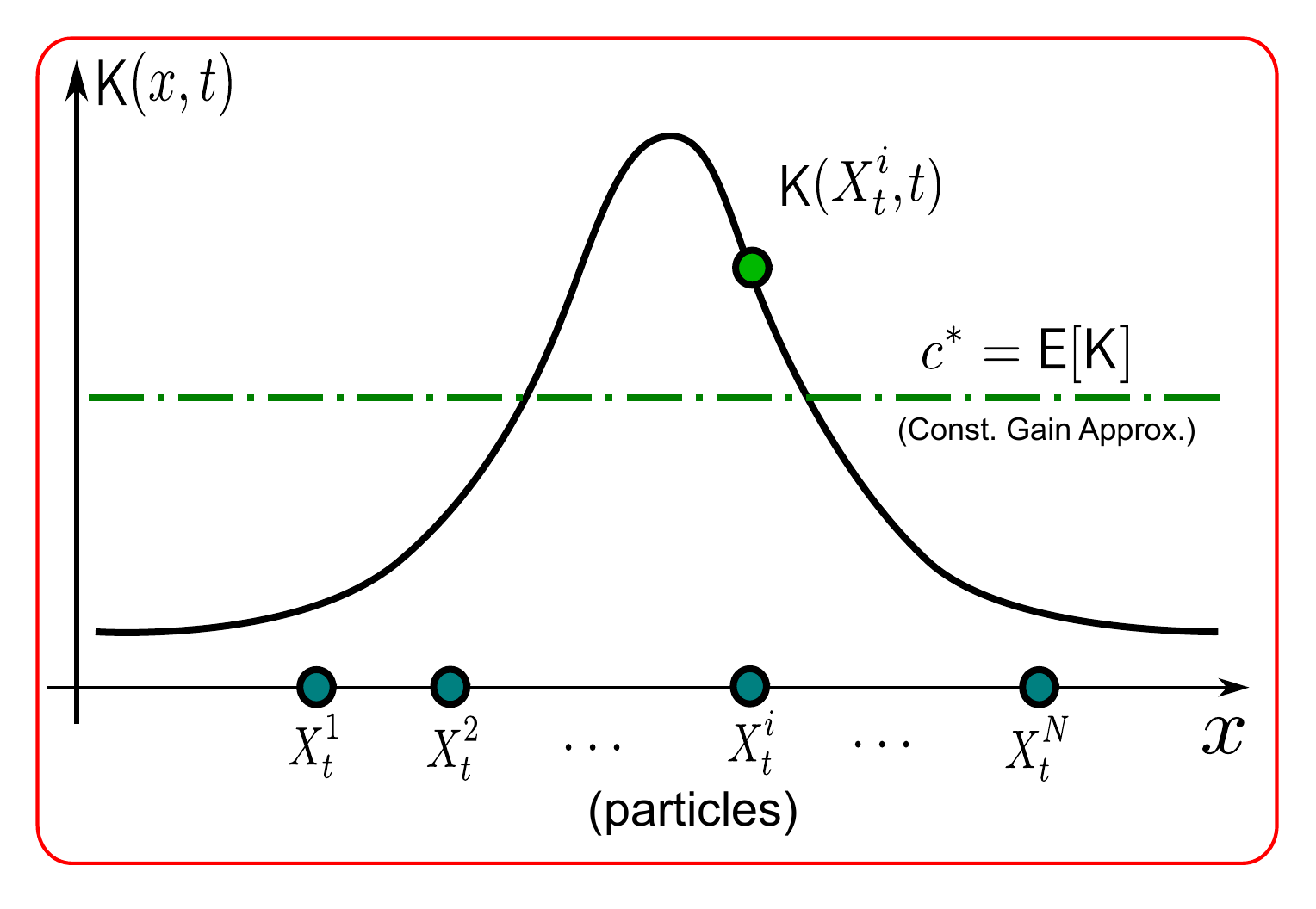}
\caption{Approximating nonlinear $\v$ by its expected value $\E[\v]$. For simplicity, the scalar case is depicted (i.e., $X_t \in \Re$).}
\label{fig:gain_func}
\vspace{-10pt}
\end{figure}

In general, the weak solution $\nabla \phi_j(\cdot, t)$ of the
BVP~\eqref{eqn:weak_form} is some, possibly non-constant, vector-valued function of
the state (see~\Fig{fig:gain_func}). The existence-uniqueness result for the
weak solution appears in~\cite{taoyang_cdc12}; a Galerkin algorithm
for numerically approximating the solution is described next:
%described in~\cite{taoyang_cdc12,taoyang_acc13,adamtilton_acc13}. 

\def\prechi{\psi}

The function $\phi_j$ is approximated as,
\begin{equation}
    \phi_j(x,t) = \sum_{l=1}^L \kappa^l_j(t) \prechi_l(x),\nonumber
\end{equation}
where $\{\prechi_l(x) \}_{l=1}^L$ are a given set of basis functions. %  These need
% to be pre-selected.

The finite-dimensional approximation of~\eqref{eqn:weak_form} is to choose
constants $\{\kappa^l_j(t)\}_{l=1}^L$ -- for each fixed time $t$ -- such that
\begin{equation}
    \sum_{l=1}^L \kappa^l_j(t) \; \Expect[ \nabla \psi_l \cdot \nabla \prechi ] =
    \Expect [  (h_j-\hat{h}_j) \psi ],\quad\forall\,\,\psi\in S,
    \label{eqn:BVP_phi_weak_expect_fd}
\end{equation}
where
$S:=\text{span}\{\prechi_1,\prechi_2,\hdots,\prechi_L\}\subset H^1(\Re^d;p)$.

Denoting $[A]_{kl} = \Expect[ \nabla \prechi_l \cdot \nabla \prechi_k ]$, $b_j^k
= \Expect [  (h_j-\hat{h}_j) \prechi_k ]$, $b_j = (b_j^1,\hdots,b_j^L)$ and $\kappa_j=(\kappa^1_j,\hdots,\kappa^L_j)$, the finite-dimensional
approximation~\eqref{eqn:BVP_phi_weak_expect_fd} is expressed as a linear matrix
equation:
\begin{equation}
    A\kappa_j = b_j.\nonumber
%    \label{eqn:lin_eqn}
\end{equation}
In a numerical implementation, the matrix $A$ and vector $b_j$ are approximated as,
\begin{align}
    [A]_{kl} & = \Expect[ \nabla \prechi_l \cdot \nabla \prechi_k ] \approx
    \frac{1}{N} \sum_{i=1}^N \nabla \prechi_l (X^i_t) \cdot
   \nabla \prechi_k (X^i_t), 
   \nonumber\\
   %\label{formula:A_ml_sample} \\
    b_j^k & = \Expect [  (h_j-\hat{h}_j) \psi_k ]  \approx \frac{1}{N} \sum_{i=1}^N
    (h_j(X^i_t) - \hat{h}_j) \psi_k(X^i_t), 
   \nonumber
   %\label{formula:b_m_sample}
\end{align}
where recall $\hat{h}_j \approx \frac{1}{N} \sum_{i'=1}^N h_j(X^{i'}_t)$.
The important point to note is that the gain function is expressed in terms of
averages taken over the population.

The constant gain approximation  is obtained by using
the coordinate functions $(x_1,x_2,\hdots,x_d)$ as basis functions.
In this case, 
\begin{equation}
\kappa_j = \E[(h_j-\hat{h}_j) x] \approx \frac{1}{N}\sum_{i=1}^N X_t^i \left(h_j(X_t^{i}) - \frac{1}{N}\sum_{i'=1}^N h_j(X_t^{i'}) \right) =: c_j^{(N)}. \label{eqn:prelim_const_gain}
\end{equation} 
Denoting $C:=[c_1^{(N)},\hdots,c_s^{(N)}]$, where $c_j^{(N)}$ is a column vector for $j\in\{1,\hdots,s\}$, the gain function is succinctly expressed as:
\begin{equation}
\v = C.\label{eqn:const_gain_prelim}
\end{equation}
We refer to this solution as the {\em {constant gain approximation}}.

\begin{remark}
There is also a variational interpretation of these solutions.  The
constant gain approximation, formula~\eqref{eqn:const_gain_prelim}, is
the best -- in the least-square sense -- constant approximation of the
gain function (see~\Fig{fig:gain_func}).  Precisely, consider the following least-square
optimization problem:
\begin{equation*}
c_j^\ast = \arg \min_{c_j \in \Re^d} \E[|\nabla \phi_j - c_j|^2].
\end{equation*}
By using a standard sum of square argument, we have
\[
c_j^\ast = \E[\nabla \phi_j].
\]
Even though $\phi_j$ is unknown, a closed-form formula for constant
vector $c_j^\ast$ can easily be obtained by using~\eqref{eqn:weak_form}. Specifically, by substituting $\psi(x) = x=(x_1,x_2,\hdots,x_d)$ in~\eqref{eqn:weak_form}:
\begin{align}
\E[\nabla \phi_j] = \E[(h_j-\hat{h}_j)\psi] = \int_{\Re^d} (h_j(x) - \hat{h}_j)\; x\;p(x,t)\ud x.\nonumber
\end{align}
This is also the first equality in~\eqref{eqn:prelim_const_gain}.

Likewise, the Galerkin solution is the optimal least-square
approximation in the function space $S$. \qed
\end{remark}

\begin{remark}
It is noted that if $p$ is Gaussian and $h$ is linear then, in the
limit as $N\rightarrow\infty$, the constant gain approximation equals the
Kalman gain (see Sec.~\ref{subsec:comp_pda}).  \qed
\end{remark}

\section{Feedback Particle Filter with Data Association Uncertainty}
\label{sec:PDAFPF}

In this section, we describe the probabilistic data association-feedback
particle filter (PDA-FPF) for the problem of tracking a single target with multiple observations. The filter for multiple targets is obtained as an extension, and described in
\Sec{sec:JPDAFPF}.

\subsection{Problem statement, assumptions and notation}
\label{sec:prob}

% At time $t$,
% we assume $M$ targets and a random number $m_t$ of observations.
% Some of the observations may be due to clutter.

The following notation is adopted:
\begin{romannum}
  \item At time $t$, the target state is denoted by
      $X_t\in\Re^d$.
  \item At time $t$, the observation vector $
      \underline{Z}_t := ( {Z}_t^1, {Z}_t^2, \hdots,
      {Z}_t^{\M})$, where $\M$ is assumed fixed and
      $Z_t^m\in\Re^s$ for $m \in \{1,\hdots, M\}$.
  \item At time $t$, the association random variable
      is denoted as $A_t\in \{0,1,\hdots,\M\}$.  It is used
      to associate one observation to the target: $A_t=m$
      signifies that the $m^{\text{th}}$-observation
      $Z_t^m$ is associated with the target, and $A_t=0$
      means that all observations at time $t$ are due to clutter. It is assumed that the target can give rise to at most one detection. We set the gating and detection probability to be $1$ for the ease of presentation. 
\end{romannum}

The following models are assumed for the three stochastic processes:
\begin{romannum}
\item The state $X_t$ evolves according to a nonlinear
    stochastic differential equation (SDE) of the form~\eqref{eqn:Signal_Process}:
\begin{equation}
\ud X_t = a(X_t)\ud t + \ud B_t,
\label{eqn:Signal_Process_Target}
\end{equation}
where the initial condition $X_0$ is drawn from a known prior distribution $p^*(x,0)$.
\item The association random process $A_t$ evolves as
    a jump Markov process in continuous-time:
\begin{equation}
{\sf P}\{A_{t+\Delta t}=m'|A_{t}=m\} = q \Delta t +
o(\Delta t),\quad m'\ne m,\label{eqn:mc_q}
\end{equation}
%The initial distribution ${\sf P}\{[A_0=m]\}=\frac{1}{\M+1}$.
where, for the ease of presentation, the transition rate is assumed to
be a constant $q$. The initial distribution is denoted as $\beta_0$.
It is assumed to be uniform.
\item $A_t$ and $X_t$ are assumed to be mutually
    independent.
\item At time $t$, the observation model is given by, for $m = 1,\hdots,M$:
\begin{equation}
\ud Z_t^{m} = 1_{[A_t=m]}\, h(X_t) \ud t + \ud W_t^m,\label{eqn:observ_model_seperate}
\end{equation}
where $\{W_t^m\}_{m=1}^M$ are mutually
independent standard Wiener processes and
\begin{equation*}
1_{[A_t=m]} :=\begin{dcases*}
        1  &  if $A_t=m$ \\
        0 & otherwise.
        \end{dcases*}
\end{equation*}
\end{romannum}

The problem is to obtain the posterior distribution of
${X}_t$ given the history of observations
(filtration) $\mathcal{\underline{Z}}_t :=
\sigma(\underline{Z}_{\tau}:\tau \le t)$. 

\begin{remark}
The equivalent ODE model to~\eqref{eqn:Signal_Process_Target} and \eqref{eqn:observ_model_seperate} is:
\begin{align}
\frac{\ud X_t}{\ud t} &= a(X_t) + \dot{B}_t,\nonumber\\
Y_t^m &= 1_{[A_t=m]}\, h(X_t)  + \dot{W}_t^m,\nonumber
\end{align}
where $\dot{B}_t$, $\{\dot{W}_t^m\}$ are independent white noise processes and $Y_t^m \doteq \frac{\ud Z_t^m}{\ud t}$.
\qed
\end{remark}

\begin{remark}
There are two differences between the clutter model assumed here and
related models used in standard discrete-time PDA algorithm~\cite{Bar-Shalom_IEEE_CSM,Bar-Shalom_book_88}:
\begin{romannum}
\item In this paper, clutter observations are modeled as a white-noise
  process in the whole space.  In the discrete-time literature, the
  standard model assumes   
  clutter observations to arise from a uniform
  distribution in a certain ``coverage area''. Related Gaussian models for
  clutter have also been considered, e.g.,~\cite{ConteLopsRicci95,WangNehorai06}.       
\item In this paper, the number of observations $M$ is assumed fixed.  
\end{romannum}
The fixed number of observations assumption helps simplify the
presentation and notation in 
the continuous-time setting of this paper.  The core algorithms can be
generalized in a straightforward manner to handle the varying number
of observations.  Concerning the Gaussian clutter model, we will
provide comparisons in the following section, in the numerical example
in Sec.~\ref{subsec:stt} and in
Appendix~\ref{apdx:discretized_assoc_filter}. 
\qed
 
% For the derivation of PDA-FPF in this paper, we model the clutter
% observations as white noise process distributed in the whole space
% rather than a uniform distribution in the``validated area'' as assumed
% in typical discrete-time PDA algorithms~\cite{Bar-Shalom_book_88},
% \cite{Bar-Shalom_IEEE_CSM}. However, the choice of the clutter
% observation density can be easily replaced in the discrete-time
% implementation of PDA-FPF. More discussion appears in
% Appendix~\ref{apdx:discretized_assoc_filter}.

\end{remark}
The PDA-FPF methodology comprises of the following two parts:
\begin{romannum}
\item Evaluation of the association probability, and
\item Integration of the association probability into the
    feedback particle filter.
\end{romannum}
% The evaluation of association probability depends upon the details of
% the clutter model -- this paper considers only the Gaussian clutter
% case.  Once the association probability is computed, the filtering
% equations for the state process -- the PDA-FPF filter -- 
% do {\em not} depend upon the details of the clutter model.    

\subsection{Association probability for a single target}

The association probability is defined as the probability of
the association $[A_t=m]$ conditioned on
$\mathcal{\underline{Z}}_t$:
\begin{equation}
\beta_t^m \triangleq \Prob \{[A_t=m]|\mathcal{\underline{Z}}_t\},\quad m = 0,1,...,\M. \nonumber
%\label{eqn:def_betaj}
\end{equation}
Since the events are mutually exclusive and exhaustive,
$\sum_{m=0}^{\M} \beta_t^m = 1$.

For the single-target-multiple-observation model described
above, the filter for computing association probability is
derived in Appendix~\ref{apdx:association_filter_pda}.
It is of the following form: For $m \in \{1,\hdots,M\}$,
% \begin{align}
% \left\{
% \begin{array}{l}
% \ud \beta_t^m  = q \left[1 - (M+1) \beta_t^m\right] \ud t+ \beta_t^m \hat{h}^T  \sum_{j=1}^{\M} \beta_t^j(\ud Z_t^m - \ud Z_t^j) + |\hat{h}|^2 \beta_t^m \sum_{j=1}^{\M} \beta_t^j (\beta_t^j - \beta_t^m)\ud t, \quad  m=1%,\hdots,\M \\
% \ud \beta_t^0 = q \left[1 - (M+1) \beta_t^0 \right] \ud t - \beta_t^0 \hat{h}^T \sum_{m=1}^M \beta_t^m \ud Z_t^m + \beta_t^0 |\hat{h}|^2 \sum_{m=1}^M (\beta_t^m)^2 \ud t, 
% \end{array}
% \right. \label{eqn:filter_for_beta_nonlinear}
% \end{align}
\begin{equation}
\ud \beta_t^m  = q \left[1 - (M+1) \beta_t^m\right] \ud t+ \beta_t^m \hat{h}^T \left(\ud Z_t^m - \sum_{j=1}^{\M} \beta_t^j \ud Z_t^j\right) + \beta_t^m |\hat{h}|^2 \left(\sum_{j=1}^{\M} (\beta_t^j)^2 - \beta_t^m\right)\ud t,
\label{eqn:filter_for_beta_nonlinear}
\end{equation}
and for $m=0$,
\begin{equation}
\ud \beta_t^0 = q \left[1 - (M+1) \beta_t^0 \right] \ud t - \beta_t^0 \hat{h}^T \sum_{j=1}^M \beta_t^j \ud Z_t^j + \beta_t^0 |\hat{h}|^2 \sum_{j=1}^M (\beta_t^j)^2 \ud t, 
\label{eqn:filter_beta_0}
\end{equation}
where $\hat{h} =  {\sf E} [h(X_t)|\mathcal{Z}_t]$ and
$|\hat{h}|^2 = \hat{h}^T \hat{h}$. This is
approximated by using particles:
\begin{equation*}
\hat{h} \approx \frac{1}{N}\sum_{i=1}^N h(X_t^i).
\end{equation*}
In practice, one may also wish to consider approaches to reduce
filter complexity, e.g., by assigning gating regions for the
observations;  cf., Sec.~4.2.3 in~\cite{Bar-Shalom_book_88}.
\begin{remark}

The association probability
filter~\eqref{eqn:filter_for_beta_nonlinear}-\eqref{eqn:filter_beta_0} can also be
derived by considering a continuous-time limit starting from the
continuous-discrete time filter in literature~\cite{Bar-Shalom_book_88}.
This proof appears
in Appendix~\ref{apdx:discretized_assoc_filter}.  The alternate proof
is included for the following reasons:
\begin{romannum}
\item The proof helps provide a comparison with the classical PDA
  filter.  This is
  important because some of the modeling assumptions (e.g., modeling
  of association $A_t$ via a jump Markov process), at first sight,
  may appear to be different from those considered in the classical literature.
% shows that the
%     filter~\eqref{eqn:filter_for_beta_nonlinear}
%   is in fact the continuous-time counterpart of the
%   algorithm that is used to
%   obtain association probability in the classical PDA filter.  This is
%   important because some of the modeling assumptions (e.g., modeling
%   of clutter, or of association $A_t$ via a jump Markov process), at first sight, may appear
%   to be different from those considered in the classical literature.
%,
%  which focusses on discrete-time observations.
\item The proof method suggests alternate discrete-time algorithms for evaluating
  association probability in simulations and experiments, where
  observations are made at discrete sampling times.\qed
\end{romannum}
\end{remark}

In the following, we integrate association probability into the
feedback particle filter, which is used to approximate the 
evolution of the posterior distribution.  The algorithmic structure
-- evaluation of data association probability first, followed by its
subsequent inclusion in the filter for approximating posterior -- is
motivated by the following considerations:
\begin{romannum}
\item Such an algorithmic structure mirrors the structure used in the 
classical PDA filtering
literature~\cite{Bar-Shalom_IEEE_CSM,Bar-Shalom_Proc_IEEE,Bar-Shalom_book_88}.
\item The computation of association probability depends upon the
  details of the clutter model -- the present paper describes this computation
  for the Gaussian clutter case
  (see~\eqref{eqn:filter_for_beta_nonlinear}-\eqref{eqn:filter_beta_0}).
  Once the association probability is computed, the filtering
  equation for the state process does not depend upon the details of
  the clutter model.  It is thus presented separately, and can also
  be used as such.  
\item A separate treatment is also useful while considering multiple
target tracking problems.  For such problems, one can extend
algorithms for data association in a straightforward manner,
while the algorithm for posterior remains as before. Additional
details appear in Sec~\ref{sec:JPDAFPF}.
\end{romannum}

\subsection{Probabilistic data association-feedback particle filter}
\label{sec:pda_fpf_main}

Following the feedback particle filter methodology, the model for the particle filter is
given by,
\begin{equation}
\ud X^i_t = a(X^i_t) \ud t + \ud B^i_t  + \ud U^i_t , \nonumber
%\label{eqn:particle_model}
\end{equation}
where $X^i_t \in \Re^d$ is the state for the $i^{\text{th}}$
particle at time $t$, $ U^i_t$ is its control input, and
$\{B^i_t\}_{i=1}^N$ are mutually independent standard Wiener processes.
We assume the initial conditions $\{X^i_0\}_{i=1}^N$  are
i.i.d., independent of $\{B^i_t\}$, and drawn from the initial
distribution $p^*(x,0)$ of $X_0$.  Both  $\{B^i_t\}$ and
$\{X^i_0\}$ are also assumed to be independent of $X_t$, $\underline{Z}_t$.
Certain additional assumptions are made regarding admissible
forms of control input (see~\cite{taoyang_cdc12}).

Recall that there are two types of conditional distributions of
interest in our analysis:
\begin{romannum}
\item $p(x,t)$:  Defines the conditional distribution of
    $X^i_t$ given $\UZ_t$.
\item $p^*(x,t)$: Defines the conditional distribution of
    $X_t$ given $\UZ_t$.
\end{romannum}

The control problem is to choose the control input $U^i_t$ so
that $p$ approximates $p^*$, and consequently empirical
distribution of the particles approximates $p^*$ for large
number of particles.

The evolution of  $p^*(x,t)$ is described by a modified form of the Kushner-Stratonovich (K-S) equation:
\begin{equation}
\ud p^\ast = \clL^\dagger p^\ast \ud t +
\sum_{m=1}^{\M} \beta_t^m ( h-\hat{h} )^T(\ud Z_t^m - \hat{h} \ud t)p^\ast. \label{eqn:PDA_KS}
\end{equation}
where $ \hat{h} := \E[h(X_t)|\UZ_t] = \int_{\Re^d} h(x) p^*(x,t) \ud x$, and $
\clL^\dagger$ is defined in~\eqref{eqn:dagger}.  The proof
appears in Appendix~\ref{apdx:consistency}. 
%see also~\cite{taoyang_acc12}.

The main result of this section is to describe an explicit
formula for the optimal control input, and to demonstrate that
under general conditions we obtain an exact match:  $p=p^*$
under optimal control. The optimally controlled dynamics of the
$i^{\text{th}}$ particle have the following Stratonovich form,
\begin{align}
\ud X_t^i = a(X^i_t)\ud t + \ud B_t^i + \underbrace{\sum_{m=1}^{\M} \beta_t^m \, \v(X_t^i,t) \circ \ud \Inov^{i,m}_t}_{\ud U_t^i},\label{eqn:PDA_FPF}
\end{align}
where $\Inov^{i,m}_t$ is a modified form of the
$\emph{innovation process}$,
\begin{equation}
\ud \Inov^{i,m}_t := \ud Z_t^m - \left[\frac{\beta_t^m}{2} h(X_t^i) + \left(1-\frac{\beta_t^m}{2}\right)\hat{h}\right]\ud t,\label{eqn:pda-inov}
\end{equation}
where $\hat{h} := {\sf E} [h(X_t^i)|\UZ_t] = \int_{\Re^d} h(x)
p(x,t) \ud x$.  

The gain function $\v = \left[\nabla \phi_1,\hdots,\nabla \phi_s\right]$ is a solution of the
E-L BVP~\eqref{eqn:EL_phi_prelim}: For $j = 1,\hdots, s$,
\begin{equation}
\label{eqn:EL_phi_pdaf}
\nabla \cdot (p(x,t) \nabla \phi_j(x,t) )  = - (h_j(x)-\hat{h}_j) p(x,t).
\end{equation}

The evolution of $p(x,t)$ is described by the forward
Kolmogorov operator for~\eqref{eqn:PDA_FPF}: See Appendix~\ref{apdx:consistency} for
the equations.

The following theorem shows that the two evolution equations
for $p$ and $p^\ast$ are identical. The proof appears in
Appendix~\ref{apdx:consistency}.

\begin{theorem}\label{thm:kushner}
Consider the two evolutions for $p$ and $p^\ast$, defined
according to the Kolmogorov forward equation~\eqref{eqn:mod_PDA_FPK} and modified K-S
equation~\eqref{eqn:PDA_KS}, respectively. Suppose that the
gain function $\v(x,t)$ is obtained according
to~\eqref{eqn:EL_phi_pdaf}. Then provided $p(x,0) =
p^\ast(x,0)$, we have for all $t \geq 0$, $p(x,t) =
p^\ast(x,t)$. \qed
\end{theorem}

\medskip

\begin{example} 
Consider the problem of tracking a target evolving according to the dynamic model~\eqref{eqn:Signal_Process_Target}, in the presence of clutter. At each time $t$, a single observation is obtained according to the observation model:
\begin{equation}
\ud Z_t = 1_{[A_t=1]}\, h(X_t) \ud t + \ud W_t,\nonumber
\end{equation}
where the association variable is denoted as $A_t \in \{0,1\}$: $A_t = 1$ signifies the event that the observation originates from the target and $A_t = 0$ means the observation is due to clutter (model of false alarm).

Let $\beta_t$ denote the observation-to-target association probability at time $t$. The probability that the observation originates from the clutter is therefore $1-\beta_t$.

For this problem, the feedback particle filter is given by:
\begin{align}
\ud X^i_t = a(X^i_t) \ud t + \ud B^i_t \; + \; \underbrace{\beta_t \,
  \v(X^i_t,t) \circ \ud \Inov^i_t}_{\ud U_t^i},
\label{eqn:particle_filter_clutter}
\end{align}
where the innovation error $\Inov^i_t$ is given by,
\begin{equation}
\ud \Inov^i_t := \ud Z_t - \left[\frac{\beta_t}{2} h(X^i_t) +
  \left(1-\frac{\beta_t}{2}\right) \hat{h} \right] \ud t.
\label{e:in_clutter}
\end{equation}

For the two extreme values of $\beta_t$, the filter
reduces to the known form:
\begin{romannum}
\item If $\beta_t=1$, the observation is associated
    with the target with probability $1$. In this case, the
    filter is the same as
    FPF~\eqref{eqn:prelim_fpf} presented
    in Sec.~\ref{sec:prelim}.
\item If $\beta_t=0$, the observation carries no
    information and the control input $\ud U^i_t =0$.
\end{romannum}

For $\beta_t\in(0,1)$, the control is more interesting. The
remarkable fact is that the innovation error-based feedback
control structure is preserved.  The association probability
serves to modify the formula for the gain function and the
innovation error:
\begin{romannum}
\item The gain function is effectively reduced to
    $\beta_t \v(X^i_t,t)$.  That is, the control gets less
    aggressive in the presence of possible false alarms due
    to clutter.
\item The innovation error is given by a more
    general formula~\eqref{e:in_clutter}.  The optimal
    prediction of the $i^{\text{th}}$ particle is now a
    weighted average of  $h(X^i_t)$ and the population
    prediction
$\hat{h} \approx \frac{1}{N} \sum_{j=1}^N h(X^j_t)$.
  Effectively, in the presence of possible false alarms, a
  particle gives more weight to the population in computing
  its innovation error.
\end{romannum}
\end{example}

\begin{remark}
The ODE model for the PDA-FPF~\eqref{eqn:PDA_FPF}-\eqref{eqn:pda-inov} is given by:
\begin{align}
\frac{\ud X_t^i}{\ud t} = a(X^i_t) + \dot{B}_t^i + \sum_{m=1}^{\M} \beta_t^m \, \v(X_t^i,t) \cdot \left( Y_t^m - \left[\frac{\beta_t^m}{2}h(X_t^i) + \left(1 - \frac{\beta_t^m}{2}\right)\hat{h}\right] \right),\label{eqn:ode_pda_fpf}
\end{align}
where $\{\dot{B}_t^i\}$ are independent white noise processes and $Y_t^m \doteq \frac{\ud Z_t^m}{\ud t}$. The gain function $\v$ is again a solution of the E-L BVP~\eqref{eqn:EL_phi_pdaf}.
\qed
\end{remark}

Table~\ref{tbl:diff} provides a comparison of the PDA-FPF and the classical Kalman filter-based PDA filter (\cite{Bar-Shalom_book_88,Bar-Shalom_IEEE_CSM}). The ODE model is adopted for the ease of presentation. The main point to note is that the feedback particle
filter has an innovation error-based structure: In effect, the $\ith$
particle makes a prediction $\hat{Y}_t^{i,m}$ as a weighted-average
of $h(X_t^{i})$ and $\hat{h}$.  This is then used to compute an
innovation error $I_t^{i,m}$.  The Bayes' update step involves gain
feedback of the innovation error.      

{\renewcommand{\baselinestretch}{1.5}
 \begin{table*}[htp]
   %\hspace{-2 in}
    \centering
    \caption{Comparison of the nonlinear PDA-FPF algorithm with the
      linear PDAF algorithm}
    \label{tbl:diff}
    \begin{tabular}{|l||c|c|}
\hline
& PDA filter & PDA-FPF\\
\hline
\hline
State model & $\dot{X}_t = A X_t + \dot{B}_t$ & $\dot{X}_t = a(X_t)+  \dot{B}_t$\\
\hline
Observation model & $Y_t = H X_t  +  \dot{W}_t$ & $Y_t = h(X_t)+  \dot{W}_t$\\
\hline
\multicolumn{1}{|l||}{\multirow{1}{*}{Assoc. Prob.}} &
        \multicolumn{2}{l|}{$\beta_t^m := \Prob\{[A_t = m]|\mathcal{\underline{Z}}^t\}$}  \\
% \cline{2-3}
% &  & \\
\hline
Prediction   & $\hat{Y}_t =H \hat{X}_t$ & $\hat{Y}_t^{i,m} = \frac{\beta_t^m}{2}h(X_t^{i}) + (1 - \frac{\beta_t^m}{2})\hat{h}$\\
\hline
Innovation error   & $I_t^m = Y_t^m- \hat{Y}_t$ & $I_t^{i,m} = Y_t^m - \hat{Y}_t^{i,m}$ \\
\hline
Feedback control   & $U_t^m = \v_g I_t^m $ & $U_t^{i,m} = \v(X_t^{i},t) I_t^{i,m}$ \\
\hline
Gain & $\v_g$: Kalman Gain & $\v(x,t)$: Sol. of a BVP~\eqref{eqn:EL_phi_pdaf} \\
\hline
Control input & $\sum_{m=1}^\M\beta_t^m U_t^m$ & $\sum_{m=1}^\M \beta_t^m U_t^{i,m}$\\
\hline
%Formula for $\beta_t^m$ & $ \propto \exp\left[- \|\ud Z_t^m-\ud \hat{Z}_t\|^2\right]$ & $\propto \frac{1}{N}\sum_{i=1}^N \exp\left[- \|\ud Z_t^m-h(X_t^{i})\ud t\|^2\right]$\\
%\hline
Filter equation & $\dot{\hat{X}}_t = A \hat{X}_t + \sum_m \beta_t^m
U_t^m$ & $\dot{X}_t^i = a(X_t^i) +  \dot{B}_t^i + \sum_m \beta_t^m U_t^{i,m}$\\
\hline
    \end{tabular}
\end{table*}
\renewcommand{\baselinestretch}{1}
}

\subsection{Algorithm}
\label{sec:algo}

For implementation purposes, we use the ODE form of the filter
(see~\eqref{eqn:ode_pda_fpf}) together with an Euler
discretization. The resulting discrete-time algorithm appears in
Algorithm~1.  

At each discrete time-step, the algorithm requires approximation
of the gain function.  The constant gain approximation algorithm
is included in Algorithm~1.  One can readily substitute another
algorithm for approximating gain function.  

The evaluation of association probability is based on the Gaussian
clutter assumption.  For this or other types of clutter models, one
could also substitute a discrete-time algorithm for association
probability (see Remark~\ref{rem:remark-DT-DAP} in Appendix~\ref{apdx:discretized_assoc_filter}).

% In practice, one can use more numerically accurate algorithm to approximate the gain function. An algorithm based on sum-of-Gaussian approximation of density appears in our conference paper~\cite{taoyang_cdc11}. In the application example presented in Sec. V-C of~\cite{taoyang_TAC12}, the gain function is approximated by using Fourier series.

Note that the PDA-FPF algorithm propagates an ensemble of particles
approximating the posterior distribution at each time step.  From the
posterior, one can estimate any desired statistic for the state.  The
formula for the mean ($X_t^{\text{est}}$) is included in Algorithm~1.

\begin{algorithm}
\label{algorithm:PDAFPF_algo}
\caption{PDA-FPF for tracking single target in clutter:}
\begin{algorithmic}[1]
\STATE {INITIALIZATION}
\FOR{$i=1$ to $N$}
\STATE Sample $X_0^{i}$ from $p(x,0)$
\ENDFOR
\FOR{$m=0$ to $M$}
\STATE Set $\beta_0^m = \frac{1}{M+1}$.
\ENDFOR
\STATE $p^{(N)}(x,0) = \frac{1}{N} \sum_{i=1}^N \delta_{X_0^{i}}(x)$
\STATE $X_{0}^{\text{est}} = \frac{1}{N}\sum_{i=1}^N X_0^i$
\end{algorithmic}
\begin{algorithmic}[1]
\STATE {ITERATION} [$t$ to $t + \Delta t$]
\FOR{$i=1$ to $N$}
\STATE Sample a Gaussian random vector $\Delta V$
\STATE Calculate $\hat{h} = \frac{1}{N}\sum_{i=1}^N h(X_t^{i})$
\STATE Calculate $|\hat{h}|^2 = \hat{h}^T \hat{h}$
\STATE Calculate the gain function $\v = \frac{1}{N}\sum_{i=1}^N X_t^i (h(X_t^i) -\hat{h})^T$
\STATE $X_{t+\Delta t}^{i} = X_t^{i} + a(X_t^{i}) \Delta  t  + \Delta V\sqrt{\Delta t} $
\FOR{$m=1$ to $M$}
\STATE Calculate $\Delta I_t^{i,m} = \Delta Z_t^m - \left[\frac{\beta_t^m}{2}h(X_t^i) + (1 - \frac{\beta_t^m}{2}) \hat{h}\right] \Delta t$
\STATE $X_{t+\Delta t}^{i} = X_{t+\Delta t}^{i} + \beta_t^{m} \v \Delta I_t^{i,m}$
%\STATE Calculate $\Delta \beta_t^m$ using~\eqref{eqn:filter_for_beta_nonlinear}
\STATE $\Delta \beta_t^m  = q \left[1 - (M+1) \beta_t^m\right] \Delta t+ \beta_t^m \hat{h}^T \left(\Delta Z_t^m - \sum_{j=1}^{\M} \beta_t^j \Delta Z_t^j\right) + \beta_t^m |\hat{h}|^2 \left(\sum_{j=1}^{\M} (\beta_t^j)^2 - \beta_t^m\right)\Delta t$
\STATE $\beta_{t+\Delta t}^m = \beta_t^m + \Delta \beta_t^m$.
\ENDFOR
\ENDFOR
\STATE $t = t+\Delta t$
\STATE $p^{(N)}(x,t) = \frac{1}{N} \sum_{i=1}^N \delta_{X_t^{i}}(x)$
\STATE $X_{t}^{\text{est}} = \frac{1}{N}\sum_{i=1}^N X_t^i$
\end{algorithmic}
\end{algorithm}

\subsection{Example: Linear case}
\label{subsec:comp_pda}
In this section, we illustrate the PDA-FPF with the aid of a linear example. The example also serves to provide a comparison to the classic PDA filter. Consider the following linear model:
\begin{subequations}
\begin{align}
\ud X_t  &= A \;X_t \ud t +  \ud B_t,\label{eqn:dyn_lin}\\
\ud Z_t^m &= 1_{[A_t = m]} \, H \; X_t \ud t + \ud W_t^m,\label{eqn:obs_lin}
\end{align}
\end{subequations}
where $A$ is a $d\times d$ matrix and $H$ is an $s \times d$ matrix.  

The PDA-FPF is described by~\eqref{eqn:PDA_FPF}-\eqref{eqn:EL_phi_pdaf}. If we assume the initial distribution $p^\ast
(x,0)$ is Gaussian with mean $\mu_0$ and covariance matrix $\Sigma_0$, then the following lemma provides the solution of the gain function $\v(x,t)$ in the linear case.

\begin{lemma}
\label{lem:linear_pda}

Consider the linear observation~\eqref{eqn:obs_lin}. Suppose $p(x,t)$ is assumed to be
Gaussian with mean $\mu_t$ and variance $\Sigma_t$, i.e., $p(x,t)=
\frac{1}{(2 \pi)^{\frac{d}{2}} |\Sigma_t|^{\frac{1}{2}}}
\exp\left[-\frac{1}{2}(x - \mu_t)^T\Sigma_t^{-1}(x - \mu_t)\right]$. Then the solution of the
E-L BVP~\eqref{eqn:EL_phi_pdaf} is given by:
\begin{align}
\v(x,t) = \Sigma_t H^T 
\label{eqn:linsol_v}
\end{align}
\qed
\end{lemma}

The formula~\eqref{eqn:linsol_v} is verified by direct
substitution in~\eqref{eqn:EL_phi_pdaf} where the
distribution $p$ is Gaussian.

The linear PDA-FPF is then given by,
\begin{align}
\ud X^i_t = A \; X^i_t \ud t +   \ud B^i_t + \Sigma_t H^T \sum_{m=1}^{\M} \beta_t^m \left(\ud Z_t^m - H \left[\frac{\beta_t^m}{2} X^i_t + \left(1-\frac{\beta_t^m}{2}\right) \mu_t\right] \ud t\right).
\label{eqn:pda_fpf_lin}
\end{align}
Note that there is no Wong-Zakai correction term since the gain function is constant. 

The following theorem states that $p=p^*$ in this case. That is, the
conditional distributions of $X_t$ and $X_t^i$ coincide. The proof is
a straightforward extension of the proof in Appendix.~\ref{apdx:consistency}, and is thus omitted.

\begin{theorem}
\label{thm_pda_linear} Consider the single target tracking
problem with a linear model defined by the state-observation
equations (\ref{eqn:dyn_lin},\ref{eqn:obs_lin}).  The PDA-FPF
is given by~\eqref{eqn:pda_fpf_lin}. In this case the posterior
distributions of $X_t$ and $X_t^i$ coincide, whose conditional
mean and covariance are given by the following,
\begin{align}
\ud \mu_t  &= A \mu_t \ud t  + \Sigma_t H^T  \sum_{m=1}^{\M} \beta_t^m \left(\ud Z_t^m - H \mu_t \ud t\right) \label{eqn:mod1}\\
 \frac{\ud \Sigma_t}{\ud t} &= A \Sigma_t + \Sigma_t A^T + I - \sum_{m=1}^{\M} (\beta_t^m)^2\, \Sigma_t H^T H \Sigma_t . \label{eqn:mod2}
\end{align}
\qed
\end{theorem}

The filter for association probability $\beta_t^m$ in the
linear Gaussian case easily follows from
using~\eqref{eqn:filter_for_beta_nonlinear}.  It is of the
following form: for $m = 1,\hdots,M$,
\begin{align}
\ud \beta_t^m  = q \left[ 1 - (M+1) \beta_t^{m} \right]\ud t  + \beta_t^m (H\mu_t)^T\left(\ud Z_t^m - \sum_{j=1}^M \beta_t^j \ud Z_t^j\right) + \beta_t^m (H \mu_t)^T (H\mu_t)\left(\sum_{j=1}^M (\beta_t^j)^2 -\beta_t^m\right)\ud t.
\label{eqn:filter_for_beta_linear}
\end{align}

In practice $\{\mu_t, \Sigma_t\}$ in \eqref{eqn:pda_fpf_lin}-\eqref{eqn:filter_for_beta_linear} are approximated as sample
means and sample covariances from the ensemble
$\{X^i_t\}_{i=1}^N$.
\begin{equation}
\begin{aligned}
\mu_t & \approx \mu_t^{(N)} := \frac{1}{N} \sum_{i=1}^N X^i_t,\\
\Sigma_t & \approx \Sigma_t^{(N)} := \frac{1}{N-1} \sum_{i=1}^N \left(X^i_t - \mu_t^{(N)}\right)\left(X^i_t - \mu_t^{(N)}\right)^T.
\end{aligned}
\nonumber
%\label{e:mut_sigmat_approx}
\end{equation}

\begin{remark}
The ODE model for the PDA-FPF~\eqref{eqn:pda_fpf_lin} is given by:
\begin{align}
\frac{\ud X_t^i}{\ud t} &= A \; X^i_t  +   \dot{B}^i_t + \Sigma_t H^T \sum_{m=1}^{\M} \beta_t^m \cdot \left(Y_t^m - H \left[\frac{\beta_t^m}{2} X^i_t + \left(1-\frac{\beta_t^m}{2}\right) \mu_t\right] \right).\nonumber
\end{align}
where $\{\dot{B}_t^i\}$ are independent white noise processes and $Y_t^m \doteq \frac{\ud Z_t^m}{\ud t}$. 
\qed
\end{remark}

%%%%%%%%%%%%%%%%%%%%%%

\section{Multiple Target Tracking using Feedback Particle Filter}
\label{sec:JPDAFPF}

In this section, we extend the PDA-FPF to the multiple target tracking
(MTT) problem. The resulting filter is referred to as the {\emph {joint
    probabilistic data association-feedback particle filter
    (JPDA-FPF)}}. For notational ease, we assume a clutter free
scenario where all observations originate from the targets (i.e.,
there are $M$ targets and $M$ observations at each time). The clutter can be handled in a manner similar to the basic PDA-FPF.

The MTT problem is introduced
in~\Sec{sec:JPDA_formulation}.  The algorithm is described {\em only}
for the special case of two targets with two observations ($M=2$),
in~\Sec{sec:JPDA_example}.  The algorithm for the general case is
conceptually analogous but notationally
cumbersome.  It appears in
Appendix.~\ref{apdx:jpda_algo}.

\subsection{Problem statement, assumptions and notation}
\label{sec:JPDA_formulation}
The following notation is adopted:

\begin{romannum}
\item There are $\N$ distinct targets. The set of targets is denoted
  by the 
    index set $\Nscr = \{1,2,\hdots, \N\}$. The set of permutations of $\Nscr$
    is denoted by $\Pscr(\Nscr)$, whose cardinality $|\Pscr(\Nscr)| = \N!$. A typical element of $\Pscr(\Nscr)$ is denoted as
    $\ua = (\alpha^1, \alpha^2, \hdots, \alpha^{\N})$. 
%where the index $\alpha^m \in \Nscr$ for $m \in \Nscr$ and $\alpha^m \neq \alpha^{m'}$ for $m\neq m'$.

\item At time $t$, the state of the $\nth$ target is denoted as $X_t^n\in\R^d$
    for $n \in \Nscr$. 

\item At time $t$, there is exactly one observation per target for a total of
    $\N$ observations, all of which are available in a centralized fashion. % The case of multiple independent sensors can also be easily handled, as will be illustrated by a numerical example in~\Sec{sec:multi_sensor}.
    The observation vector is denoted by $ \underline{Z}_t :=
    ({Z}_t^{1}, {Z}_t^{2}, \hdots,Z_t^{\N})$, where the $\mth$ entry,
    $Z_t^{m} \in\R^s$, originates from one of the targets in $\Nscr$. 

\item At time $t$, the association random vector is denoted as $\uA_t \in
    \Pscr(\Nscr)$. It is used to associate targets with the observations in a joint manner: $\uA_t
    := (\alpha_t^1, \alpha_t^2,\hdots, \alpha_t^\N) \in \Pscr(\Nscr)$ signifies
    that the observation $Z_t^1$ originates from target $\alpha_t^1$, $Z_t^2$
    originates from target $\alpha_t^2$, \ldots, $Z_t^\N$ originates from target
    $\alpha_t^\N$. 
    
\end{romannum}

\medskip

The following models are assumed for the three stochastic processes.  
\begin{romannum}
\item The dynamics of the $\nth$ target evolves according to the nonlinear SDE:
    \begin{align}
	\ud X^{n}_t &= a^{n}(X_t^n) \ud t  + 	\ud {B}^{n}_{t}\, , \label{eqn:mul_tar_model}
    \end{align}
    where $\{B^{n}_{t}\in\R^d\}$ are mutually independent
    standard Wiener processes and $n\in \Mscr$.
\item The observation model is given by: 
    \begin{align}
	\bmat{\ud Z_t^1\\ \vdots \\\ud Z_t^{\N}} &= \Psi(\uA_t) \bmat{\bearing(\vX_t^1)\\
	\vdots \\ \bearing(\vX_t^{\N})} \ud t + \bmat{\ud W_t^1\\ \vdots\\ \ud W_t^{\N}}, \label{eqn:mul_obs_model}
    \end{align}
    where $\Psi(\uA_t)$ is the permutation matrix for association vector $\uA_t$, and $\{W_t^{m}\}_{m=1}^M$ are mutually independent standard Wiener processes, also assumed to be mutually independent with $\{{B}^{n}_t\}$.
   
\item The model for $\uA_t$ is similar to the model assumed in PDA-FPF (see~\Sec{sec:PDAFPF}) and is described by a continuous-time Markov chain:
    \begin{equation}
  \P\{[\uA_{t+\Delta t} = \ugamma]|[\uA_t = \ugamma']\} = q \Delta t + o(\Delta t),\quad \ugamma \neq \ugamma' \in \Pi(\Mscr),
  \label{eqn:def_At_joint}
  \end{equation}
where, for the ease of presentation, the transition rate is assumed to be a constant $q$. The initial distribution is assumed given.
\item $\uA_t$ and $\uX_t$ are assumed to be mutually independent.
\end{romannum}

The problem is to design $\N$ feedback particle filters, where the $\nth$
filter is intended to estimate the posterior distribution of the $\nth$ target
given the history of all un-associated observations (filtration) $\UZ_t :=
\sigma\left(\uZ_{\tau}: \tau \leq t\right)$. %The posterior distribution for the $\nth$ target is denoted as $p_n(x,t)$.

\begin{remark}
The equivalent ODE model to~\eqref{eqn:mul_tar_model} and \eqref{eqn:mul_obs_model} is:
\begin{align}
\frac{\ud X_t^n}{\ud t} &= a^n(X_t^n) + \dot{B}_t^n,\quad n = \{1,\hdots,M\}\nonumber\\
\bmat{Y_t^1\\ \vdots \\ Y_t^{\N}} &= \Psi(\uA_t) \bmat{\bearing(\vX_t^1)\\
	\vdots \\ \bearing(\vX_t^{\N})} + \bmat{\dot{W}_t^1\\ \vdots\\ \dot{W}_t^{\N}}, \nonumber
\end{align}
where $\{\dot{B}_t^n\}$, $\{\dot{W}_t^m\}$ are independent white noise processes and $Y_t^m \doteq \frac{\ud Z_t^m}{\ud t}$.
\qed
\end{remark}

The joint association probability is defined as the probability of the joint association event $[\uA_t = \ugamma]$ conditioned on $\UZ_t$:
\begin{equation}
\pi_t^{\ugamma} := \P\{[\uA_t = \ugamma]|\UZ_t\}, \quad \ugamma \in \Pi(\Mscr).
\label{eqn:def_joint_prob}
\end{equation}

The filter for the joint association probability $\pi_t^{\ugamma}$ is
a notationally tedious but straightforward extension of the
association filter in PDA-FPF~\eqref{eqn:filter_for_beta_nonlinear}.  It appears in Appendix~\ref{apdx:jpda_algo}.

The joint association probabilities are then used to obtain
association probability for each individual target. Once the
observation-to-target association probability is known, its
integration into the feedback particle filter is identical to 
the PDA-FPF~\eqref{eqn:PDA_FPF}. The complete set of equations for the
JPDA-FPF algorithm appear in Appendix~\ref{apdx:jpda_algo}.  The $M=2$ case is
illustrated next.
%In the following, a two-target two-observation example is discussed to illustrate the use of JPDA-FPF.

\subsection{Example: Two-target two-observation problem}
\label{sec:JPDA_example}

% In this section, a two-target two-observation problem is considered to illustrate the JPDA-FPF algorithm. 

At time $t$, the target state is denoted as
    $\underline{X}_t:=(X_t^1,X_t^2)$, and the observation vector $
    \underline{Z}_t := ({Z}_t^1, {Z}_t^2)$. The association random vector
    is denoted as $\uA_t\in \Pscr(\{1,2\}) := \{\ugamma_1,\ugamma_2\}$, where $\ugamma_1 := (1,2)$ and $\ugamma_2 := (2,1)$. It is used to
    associate observations to targets: 
    $\uA_t=\ugamma_1$ signifies that $Z_t^1$ originates from 
    target $1$ and $Z_t^2$ from target $2$; 
    $\uA_t=\ugamma_2$ accounts for the complementary case.

The following models are assumed for the three stochastic
processes:
\begin{romannum}
\item Each element of the state vector
    $\underline{X}_t$ evolves according to a nonlinear SDE of the
    form~\eqref{eqn:Signal_Process}:
\begin{equation}
\ud X_t^n = a^n(X_t^n)\ud t + \ud B_t^n,\quad n\in\{1,2\},
\label{eqn:Signal_Process_Two_Target}
\end{equation}
where $\{B_t^1\}$,$\{B_t^2\}$ are mutually independent standard
Wiener processes.
\item The association random process $A_t$ evolves as
    a jump Markov process in continuous-time:
\begin{equation}
{\sf P}\{[\uA_{t+\Delta t}=\ugamma_{m'}]|[\uA_{t}=\ugamma_{m}]\} = q \Delta t+
o(\Delta t),\quad m'\ne m \in \{1,2\}.
\nonumber
\end{equation}
The initial distribution is assumed given.
\item $\uA_t$ and $\underline{X}_t$ are assumed to be
    mutually independent.
\item At time $t$, the observation model is given by,
\begin{equation}
\left[
\begin{array}{ccc}
\ud Z_t^1 \\
\ud Z_t^2
\end{array} \right] = \Psi (\uA_t) \left[
                                    \begin{array}{ccc}
                                    h(X_t^1)\\
                                    h(X_t^2)
                                    \end{array}\right] \ud t +\left[
                                                         \begin{array}{ccc}
                                                         \ud W_t^1 \\
                                                         \ud W_t^2
                                                         \end{array}\right],\nonumber
                                                         %\label{eqn:Two_Target_Observ_Process}
\end{equation}
where $\{W_t^1\},\{W_t^2\}$ are mutually independent Wiener processes and $\Psi(\uA_t)$ is a function mapping $\uA_t$
to a permutation matrix:
\begin{equation}
\Psi(\ugamma_1) =
\begin{bmatrix}
        I_s & 0 \\
        0 & I_s
       \end{bmatrix},\quad \Psi(\ugamma_2) = \begin{bmatrix}
                                                0 & I_s\\
                                                I_s & 0
                                                \end{bmatrix},\nonumber
                                                %\label{eqn:permuation}
\end{equation}
where $I_s$ is the $s \times s$ identity matrix.
\end{romannum}
 
The joint association probability is defined as the probability
of the joint association $[\uA_t = \ugamma]$ conditioned on $\UZ_t$:
\begin{equation}
\pi_t^m := \Prob \{[\uA_t = \ugamma_m]|\UZ_t\},\quad m \in \{1,2\}. \nonumber
%\label{eqn:jpda_def}
\end{equation}

Based on the JPDA-FPF algorithm discussed in Appendix.~\ref{apdx:jpda_algo}, the filter for joint association probability $\pi_t^1$ is:
\begin{align}
\ud \pi_t^1 = -q (\pi_t^1 - \pi_t^2) \ud t + \pi_t^1 \pi_t^2 \hat{\tilde{h}}^T (\ud Z_t^1 - \ud Z_t^2)  - (\pi_t^1 - \pi_t^2) \pi_t^1 \pi_t^2 \left|\hat{\tilde{h}}\right|^2 \ud t, 
\nonumber
%\label{eqn:jpda_dpi_1}
\end{align}
where $\tilde{h}(X_t^1,X_t^2) := h(X_t^1) -h(X_t^2)$, $\hat{\tilde{h}} := \Expect[\tilde{h}|\UZ_t]$ and
$\left|\hat{\tilde{h}}\right|^2 = \hat{\tilde{h}}^T \hat{\tilde{h}}$. The expectations are approximated by using particles. Since the joint
events are mutually exclusive and exhaustive, we have
$\sum_{m=1}^2 \pi_t^m = 1$. Using this, we have $\pi_t^2 = 1-
\pi_t^1$ and $\ud \pi_t^2 = - \ud \pi_t^1$.

The joint association probabilities $\pi_t^1,\pi_t^2$ are used
to obtain association probability for individual
target. For example, for target $1$: $\beta_t^{1,1} = \pi_t^1$,
$\beta_t^{2,1} = \pi_t^2$, where $\beta_t^{m,n}$ denotes the
conditional probability that the $\mth$ observation originates from
the $\nth$ target. Once the association probability is known, the feedback particle filter for each target is of the
form~\eqref{eqn:PDA_FPF}. The detailed algorithm for JPDA-FPF is omitted because of its similarity to PDA-FPF.

\section{Numerics} 
\label{sec:numerics}

\subsection{Single target tracking in clutter}
\label{subsec:stt}

Consider first the problem of tracking a single target in clutter.
The target evolves according to a white-noise acceleration
model:
\begin{align}
\ud X_t &= FX_t\ud t + \sigma_B \ud B_t,\label{eqn:stt_dyn}\\
\ud Z_t &= HX_t\ud t + \sigma_W \ud W_t,\label{eqn:stt_obs}
\end{align}
where $X_t$ denotes the state vector comprising of position and
velocity coordinates at time $t$, $Z_t$ is the observation
process, $\{B_t\},\{W_t\}$ are mutually independent standard
Wiener processes. The two matrices are given by:
\begin{equation*}
F=\begin{bmatrix}
                    0 &  1 \\
                    0 &  0
  \end{bmatrix},
\quad H = \begin{bmatrix}
                    1 &  0
    \end{bmatrix}.
\end{equation*}

\begin{figure}
    \centering
    \vspace{-8pt}
    \includegraphics[scale=0.6]{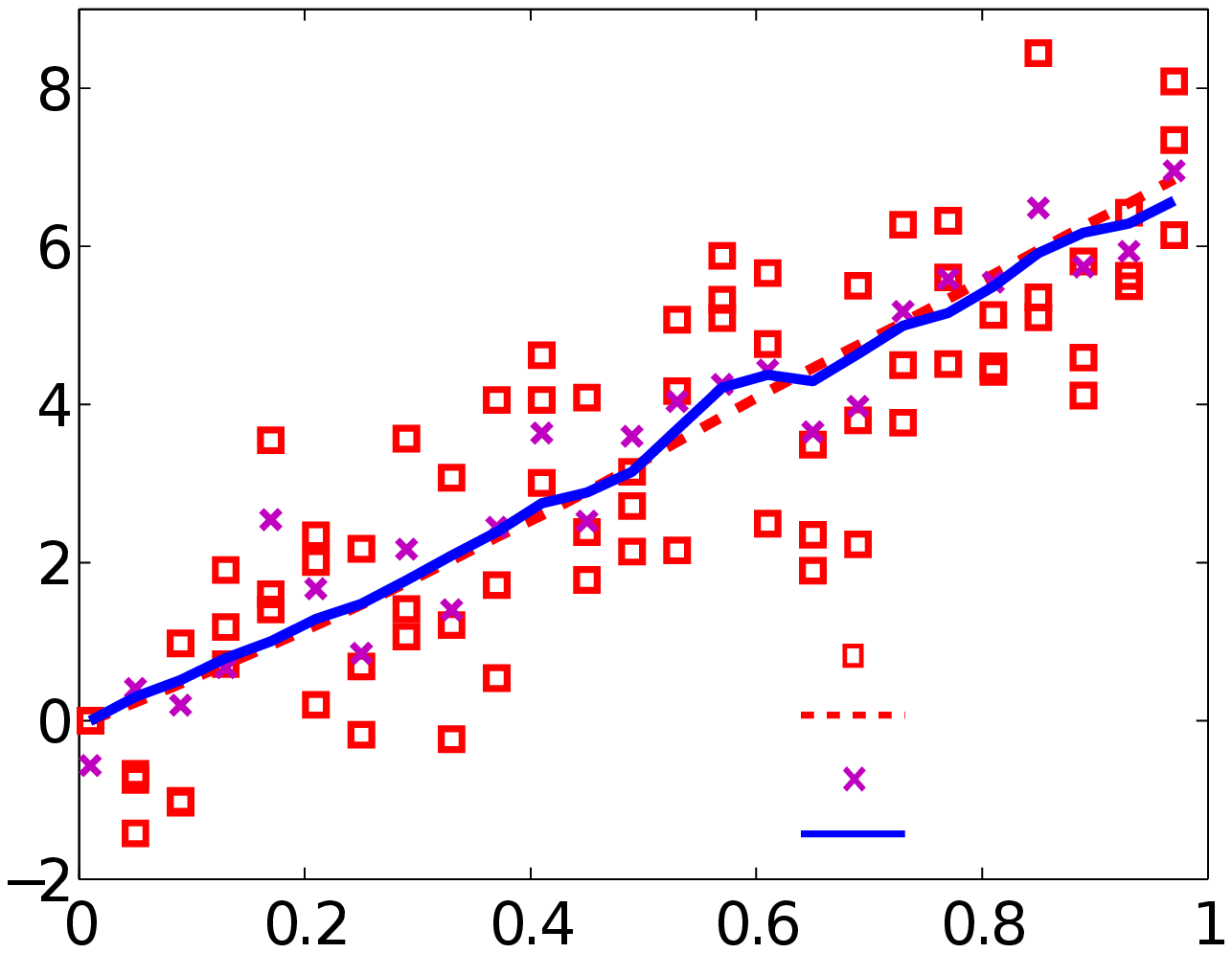}
    \caption{Simulation results of single target tracking in clutter using PDA-FPF: Comparison of estimated mean with the true trajectory.}
    \label{fig_stt}
\end{figure}

In the simulation results described next, the following
parameter values are used: $\sigma_B = (0,1)$,
$\sigma_W = 0.06$, and initial condition $X_0 =
(0,6)$.  The total simulation time 
$T=1$ and the discrete time-step $\Delta t = 0.01$. 
%The transition rate of the Markov chain~\eqref{eqn:mc_q} is $q= 10$. 
At each
discrete-time step, we assume $M=4$ observations, one due to the target and the others because of clutter. 
For the clutter observations, the standard discrete-time model
from~\cite{Bar-Shalom_IEEE_CSM} is
assumed: The clutter observations are sampled uniformly from a
coverage area of radius $2$ centered around the target.  
%The associations are not apriori known.

Figure~\ref{fig_stt} depicts the result of a single simulation: 
True target trajectory is depicted as a dashed line. At each
discrete time step, target-oriented observations are depicted
as crosses while clutter observations are depicted as squares.

The PDA-FPF is implemented according to Algorithm~1 in
Sec.~\ref{sec:algo}.  For the evaluation of association probability,
the transition rate parameter of the Markov chain~\eqref{eqn:mc_q} is $q=10$.  The output
of the algorithm -- the estimated mean trajectory -- is depicted as a
solid line. 
For the filter simulation, 
$N=1000$ particles are used.  The initial ensemble of particles are drawn from a Gaussian distribution whose mean is $\mu_0 = X_0$ and covariance matrix is $\Sigma_0 = \begin{bmatrix}
                    0.1 &  0 \\
                    0 &  0.05
  \end{bmatrix}$.

% \begin{remark}
% In spite of the Gaussian clutter model assumption in~\Sec{sec:PDAFPF}-A, we use a uniformly distributed clutter model in the discrete-time implementation of this example. The coverage area $V$ is set to be $V=4$ centered around the true observation. The assumption of Gaussian clutter is for theoretical derivation of the continuous-time data association filter~\eqref{eqn:filter_for_beta_nonlinear}-\eqref{eqn:filter_beta_0}. More detail appears in~\Appendix{apdx:discretized_assoc_filter}. \qed
% \end{remark}

\subsection{Two-target tracking problem}
\label{sec:multi_sensor}

\begin{figure*}
    \centering
    \Ebox{1}{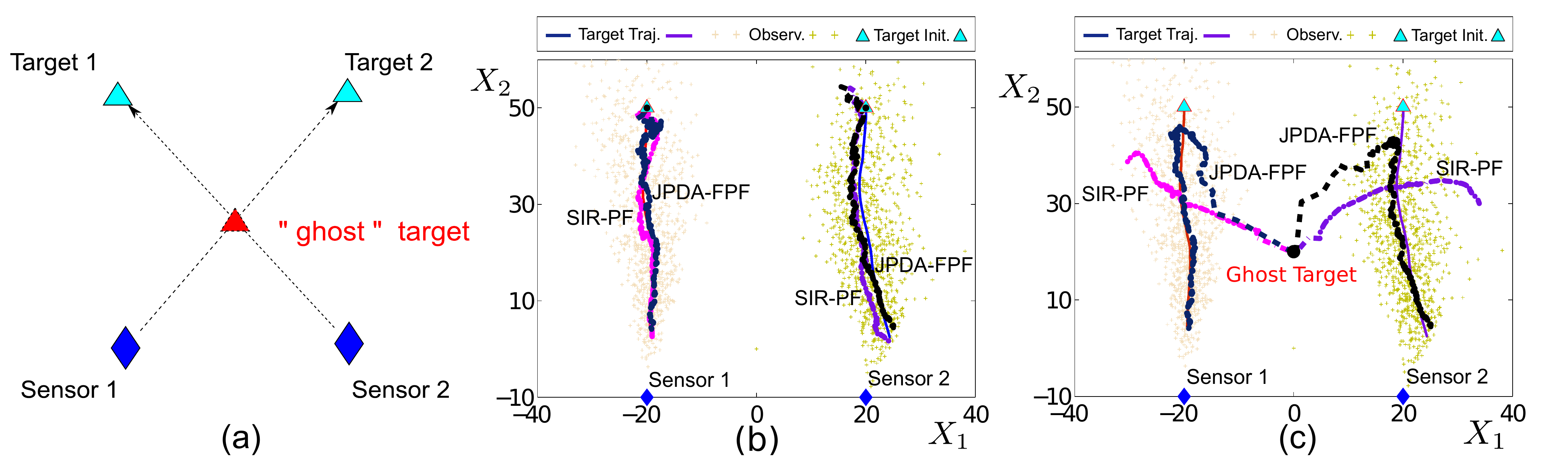}
    \caption{(a) Illustration of ``ghost'' target in the two-sensor
      two-target case: The ghost appears because of incorrect data
      association. Simulation results for (b) scenario~(i) and (c) scenario~(ii).}
    \vspace{-5pt}
    \label{fig:S_incoming}
    \vspace{-10pt}
\end{figure*}

Consider the tracking problem for two targets with two bearing-only
sensors as depicted in~\Fig{fig:S_incoming}(a).  Each target moves in a two-dimensional (2d) plane according to the standard white-noise acceleration model:
\begin{equation}
    \ud X_t = A X_t \ud t  + \Gamma \ud B_t,\nonumber
%    \label{eqn:white_noise_model}
\end{equation}
where $X :=(X_1,V_1,X_2,V_2)\in \R^4$, $(X_1,X_2)$ denotes the 2d-position of a target and $(V_1,V_2)$ denotes its velocity.  The matrices, 
\begin{equation}
A = \begin{bmatrix}
	0 & 1 & 0 & 0\\
	0 & 0 & 0 & 0\\
	0 & 0 & 0 & 1\\
	0 & 0 & 0 & 0
    \end{bmatrix},\qquad 
\Gamma = \sigma_B \begin{bmatrix}
		    0 & 0\\
		    1 & 0\\
		    0 & 0\\
		    0 & 1
		  \end{bmatrix},\nonumber
\end{equation}   
and $\{ B_t\}$ is a standard 2d Wiener process.

%\begin{figure*}
%    \centering
%    \Ebox{1}{ghost_tgt_combine_2}
%    \caption{(a) Illustration of ``ghost'' target in the two-sensor two-target case: The ghost appears because of incorrect data association. (b) Simulation results for good track initialization. (c) Simulation results for bad track initialization.}
%    \vspace{-5pt}
%    \label{fig:S_incoming}
%\end{figure*}

The two bearing-only sensors are also depicted in the figure. At time
$t$, two angle-only observations are made for each target according to the observation model: 
\begin{equation}
    \ud Z_t = h(X_t) \ud t+ \sigma_W \ud W_t,\nonumber
%    \label{eqn:num_obs_model}
\end{equation}
where $\{W_t\}$ is a standard 2d Wiener process, $h = (h_1,h_2)$ and
\begin{equation}
    h_j(x_1,v_1,x_2,v_2) = \arctan \left(\frac{x_2 - x_2^{(\text{sen}\;
    j)}}{x_1 - x_1^{(\text{sen}\; j)}}\right),\quad j \in \{1,2\},\nonumber
\end{equation}
where $(x_1^{(\text{sen}\; j)},x_2^{(\text{sen}\; j)})$ denote the position of
sensor $j$.

There exists data association uncertainty, in the sense that one cannot assign
observations to individual targets in an apriori manner. In this
particular example, faulty data association can lead to appearance of  a ``ghost''
target (see~\Fig{fig:S_incoming}(a)). The ``ghost'' target position
is identified by using the process of triangulation based on angle only observations with the two targets.

In the simulation study described next, we compare performance of the JPDA-FPF
algorithm with the sequential importance resampling particle filter
(SIR-PF) algorithm.  In the SIR-PF implementation, the data association
uncertainty is taken into account by updating the particle weights
according to the algorithm described in~\cite{Karlsson01}.

The performance comparisons are carried out for two distinct filtering
scenarios:
\begin{romannum}
\item The target tracks for the filters are initialized at the true
  target positions (depicted as black circles in~\Fig{fig:S_incoming}(b)); 
\item The target tracks for the filters are initialized at the ghost
  target position (depicted as black circles in~\Fig{fig:S_incoming}(c)). 
\end{romannum}
In the filter implementation,  particles are initialized by drawing
from a Gaussian distribution whose mean is the initial track
position and the covariance matrix $\Sigma_0 =
\text{diag}(\{10,1,10,1\})$.

Figure~\ref{fig:S_incoming} parts~(b) and (c) depict the simulation
results for the scenarios~(i) and~(ii), respectively.  In each of the
figure, the true target track is depicted together with the estimate
(mean) from the two filter implementations.  In the background, the
ensemble of observations are shown.  Each point in the ensemble is
obtained by using the process of triangulation based on two (noisy)
angle observations.

The JPDA-FPF algorithm is based on the methodology described
in~\Sec{sec:JPDA_example}, where the PDA-FPF is implemented for each
target according to Algorithm~1 in Sec.~\ref{sec:algo}.  At each time $t$, four association probabilities are approximated for each sensor:
$\ProbAsso_t^{1,1;1}$, $\ProbAsso_t^{2,1;1}$ are the two association
probabilities for target $1$,  and $\ProbAsso_t^{1,2;1}$,
$\ProbAsso_t^{2,2;1}$ are the probabilities for target 2, where
$\beta_t^{m,n;r}$ denotes the conditional probability that the $\mth$
observation of the $\rth$ sensor originates from the $\nth$ target.

% Two different initial tracks are considered: the two tracks are initialized (i). at the correct location of the target (black circles in~\Fig{fig:S_incoming}(b)); (ii). at the location of the ``ghost'' target (black circles
% in~\Fig{fig:S_incoming}(c)). Particles are initialized by drawing from a Gaussian distribution whose mean is set to be the initial track position and a covariance matrix $\Sigma_0$. Figure~\ref{fig:S_incoming}(b)-(c) depicts a sample path obtained in the Monte Carlo numerical experiments. The sensor and target locations are depicted together with a track estimate (particle means) that is approximated using the JPDA-FPF and the SIR PF algorithm. The background depicts the ensemble of observations that were
% made over the simulation run.  Each point in the ensemble is obtained by using
% the process of triangulation based on two (noisy) angle observations. The JPDA-FPF algorithm described in~\Sec{sec:JPDA_example} is implemented with the gain function approximated by~\eqref{eqn:prelim_const_gain}. At each time $t$, four association probabilities are approximated for each sensor:
% $\ProbAsso_t^{1,1;1}$, $\ProbAsso_t^{2,1;1}$ are the two association probabilities for target $1$,  and $\ProbAsso_t^{1,2;1}$, $\ProbAsso_t^{2,2;1}$ are the probabilities for target 2, where $\beta_t^{m,n;r}$ defines the conditional probability that the $\mth$ observation of the $\rth$ sensor originates from the $\nth$ target.

The simulation parameters are as follows: The two targets start at position $(-20, 50)$ and $(20, 50)$,
respectively. The ghost target location is at $(0,20)$. The initial velocity vectors are $V_0^1 = V_0^2 =
(0.0, -5.0)$.  The noise parameters $\sigma_B = 0.5$, $\sigma_W =
0.01$ and the number of particles $N = 200$. 
% The covariance matrix for
% initial particle ensemble is $\Sigma_0 =
% \text{diag}(\{10,1,10,1\})$. 
The total simulation time is $T=10$ and
the fixed discrete time-step $\Delta t = 0.01$. The transition rate of the
Markov chain~\eqref{eqn:def_At_joint} is $q = 10$. 

The comparison for the two scenarios show that: 
\begin{romannum}
\item If the filter tracks are initialized at the correct locations of
  the targets, both JPDA-FPF and SIR-PF algorithms are able to track
  the targets; 
\item If the tracks are initialized at the incorrect ghost target
  location, the SIR-PF diverges while the JPDA-FPF can adequately
  track both targets after a short period of transients.    
\end{romannum}

%\begin{figure*}
%    \centering
%    \Ebox{1}{ghost_tgt_combine_4}
%    \caption{(a) Illustration of ``ghost'' target in the two-sensor
%      two-target case: The ghost appears because of incorrect data
%      association. Simulation results for (b) scenario~(i) and (c) scenario~(ii).}
%    \vspace{-5pt}
%    \label{fig:S_incoming}
%    \vspace{-10pt}
%\end{figure*}

\begin{figure*}
    \centering
    \Ebox{0.8}{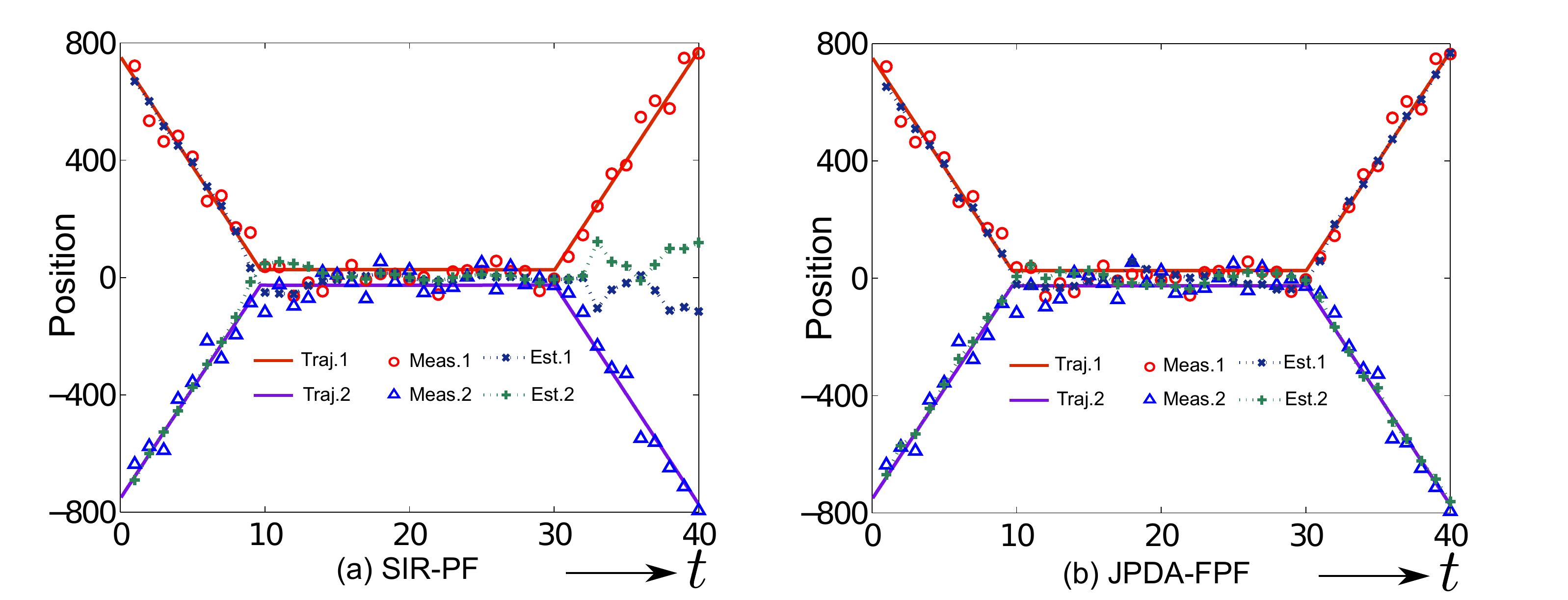}
    \caption{(a) Simulation results for SIR-PF; (b) Simulation results for JPDF-FPF.}
    \label{fig:SIR_JPDAFPF}
\end{figure*}

%-----------------------------------------------------------------------------
\subsection{Track coalescence}
%-----------------------------------------------------------------------------
Track coalescence is a well-known problem in multiple target tracking
applications. Track coalescence can occur when two closely
spaced targets move with approximately the same velocity over a period
of time~\cite{Blackman_book,BlomBloem2000}.  With a standard
implementation of the SIR-PF algorithm,
the target tracks are known to coalesce even after the targets have
moved apart~\cite{BlomBloem2000}. In the following example, simulation
results with JPDA-FPF are described for a model problem
scenario taken from~\cite{BlomBloem2006}. The results are compared
with the SIR-PF algorithm taken from~\cite{Avitzour95,Karlsson01}.

The simulation study involves two targets moving on a one-dimension
line~\cite{BlomBloem2006}. At time $t=0$, the two targets start at position
$750$ and $-750$, respectively. They move towards each other according
to a white noise acceleration model~\eqref{eqn:stt_dyn} with initial
velocity vector $-75$ and $75$, respectively.  At time $t=10$, the
distance between the two targets is smaller than $d = 50$, and they
come to a stop and maintain position for the next $20$ seconds. At
time $t=30$,  the two targets move away from each other with velocity $75$ and
$-75$, respectively. The resulting position trajectories of the two
targets (denoted as $X_t^{\text{tgt}}$) are depicted in Figure~\ref{fig:SIR_JPDAFPF}. The target
positions are observed according to the observation
model~\eqref{eqn:stt_obs}, again with data association uncertainty.

In the Monte Carlo simulation study described next, the following
noise parameters are used: $\sigma_B = (0,25)$ and $\sigma_W =10$. One
hundred Monte Carlo simulations were performed over a total simulation
period $T = 40$ with a fixed discrete time-step $\Delta t = 0.05$.  
%The transition rate of the Markov chain~\eqref{eqn:def_At_joint} is
%$q= 10$. 
At each discrete time-step, two observations are made but the
association is not known in an apriori manner.  For the filters,
$N=1000$ 
particles are initialized by drawing from a Gaussian distribution
whose mean is set to be the target's initial position and covariance
matrix $\Sigma_0 = \text{diag}(\{100,10,100,10\})$.  The JPDA-FPF is
implemented as described in~\Sec{sec:JPDA_example} and the SIR-PF
algorithm is taken from~\cite{Karlsson01}.
% In both cases, $N=1000$
% particles are used. 

\begin{figure*}
    \centering
    \Ebox{0.8}{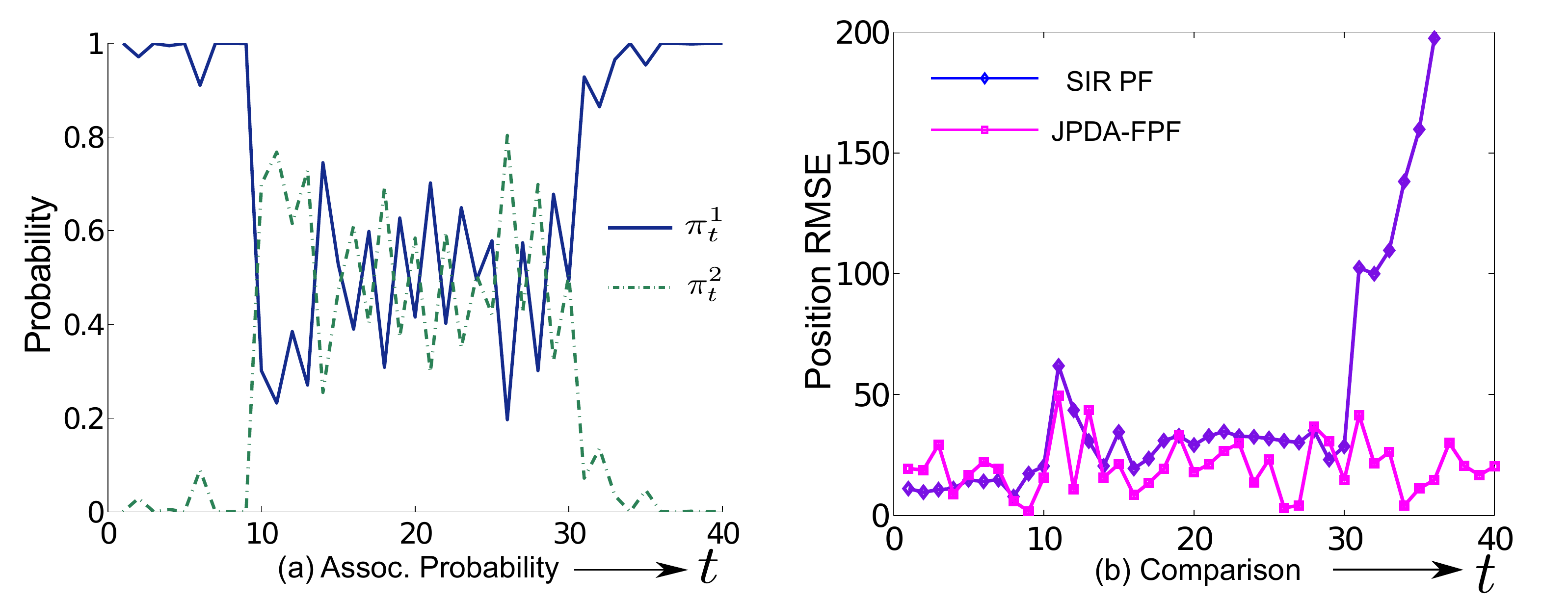}
    \caption{(a) Plot of data association probability; (b) Comparison of RMSE with SIR-PF and JPDA-FPF}
    \label{fig:Prob_RMSE}
\end{figure*}

Simulation results are summarized in
Figure~\ref{fig:SIR_JPDAFPF}-\ref{fig:Prob_RMSE}. 
Figure~\ref{fig:SIR_JPDAFPF} parts~(a) and~(b) depict the typical
results for one simulation run with SIR-PF and JPDA-FPF algorithms,
respectively.  Track coalescence is observed with the SIR-PF algorithm
while the JPDA-FPF is able to track the two targets.
Figure~\ref{fig:Prob_RMSE} (a) depicts the corresponding trajectory of
the association probability with the JPDA-FPF algorithm.  As shown in the
figure, the tracks coalesce when the two targets are close. However,
the filter is able to recover the target tracks once the targets move
away. 

Figure~\ref{fig:Prob_RMSE} (b) depicts a comparison of the root mean
square error (RMSE) in position with the SIR-PF and JPDA-FPF
algorithms.  The RMSE at time $t$ is defined as:
\begin{equation}
\text{RMSE}_t = \sqrt{\E\left[\left|X_t^\text{est}-X_t^\text{tgt}\right|^2\right]},\nonumber
%\label{eqn:rmse_def}
\end{equation}
where $\E[\cdot]$ is computed by taking the average of $100$ Monte
Carlo runs. For $t<30$, the performance of the two algorithms is
comparable.  However once the two closely spaced targets depart away
from each other, the RMSE of SIR-PF soars while the JPDA-FPF remains
relatively stable.

%
%\begin{figure*}
%    \centering
%    \Ebox{0.9}{sir_jpdafpf}
%    \caption{(a) Simulation results for SIR-PF; (b) Simulation results for JPDF-FPF.}
%    \label{fig:SIR_JPDAFPF}
%\end{figure*}

{\renewcommand{\baselinestretch}{1.5}
 \begin{table*}[htp]
   %\hspace{-2 in}
    \centering
    \caption{Comparison between SIR PF and JPDA-FPF}
    \label{tbl:rmse}
    \vspace{-5 pt}
    \begin{tabular}{|l||c|c|}
\hline
Method & Avg. RMSE & Avg. \% Tracks ``OK''\\
\hline
\hline
SIR-PF & 75.40 & 78\\
\hline
JPDA-FPF & 22.86 & 100\\
\hline
    \end{tabular}
\end{table*}
\renewcommand{\baselinestretch}{1}
}

Table~\ref{tbl:rmse} provides a summary of the performance with the
two algorithms.  The metrics are average RMSE and the average
percentage of the filter tracks being ``OK.''  The average RMSE is calculated as:
\begin{equation*}
\text{avg. RMSE} = \sqrt{\frac{1}{T}\int_0^T \E\left[\left| X_t^{\text{est}} - X_t^{\text{tgt}} \right|^2\right] \ud t}.
\end{equation*}
A track is said to be ``OK'' if the average tracking error satisfies:
\begin{equation}
\sqrt{\frac{1}{T}\int_0^T \left|X_t^{\text{est}} -
    X_t^{\text{tgt}}\right|^2 \ud t} \leq 9 \sigma_W.\nonumber
\end{equation}

The metrics used here are borrowed from~\cite{Karlsson01,BlomBloem11Fusion}.

%%%%%%%%%%%%%%%%%%%%%%%%%%%%%%%%%%%%%%%%%%%%%%%%%%%%%%%%%%%%%%%%%%%%%%%%%%%%%

\section{Conclusion}

In this paper, we introduced a novel feedback control-based particle
filter algorithm for the filtering problem with data association
uncertainty. The algorithm is referred to as the probabilistic data
association-feedback particle filter (PDA-FPF).  PDA-FPF provides for
a generalization of the Kalman filter-based PDA filter to a general
class of nonlinear non-Gaussian problems.  Our proposed algorithm
inherits many of the properties that has made the Kalman filter-based
algorithms so widely applicable over the past five decades, including
innovation error and the feedback structure
(see~\Fig{fig:fig_FPF_KF}). 

The feedback structure is important on account of the issue of
robustness.  The structural aspects are also expected to be useful for design,
integration and testing of the algorithm in a larger system involving
filtering problems (e.g., navigation systems).   

Several numerical examples were included to illustrate the theoretical
results, and provide comparisons with the SIR particle filter. We refer the
reader to our paper~\cite{taoyang_cdc13} where feedback particle
filter-based algorithms for interacting multiple models are described. 
%  For the problem
% of maneuvering target tracking, we refer the reader to our
% paper~\cite{taoyang_cdc13} where the feedback particle filter-based
% algorithm are extended to nonlinear filtering problem of stochastic
% hybrid systems.

% The main purpose of this paper has been to introduce the algorithm and
% structural features of the PDA-FPF for filtering problems with data
% association uncertainty. A complete numerical comparison, in terms of
% tracking performance with the existing approaches, for MTT applications
% is planned for future publications.

%%%%%%%%%%%%%%%%%%%%%%%%%%%%%%%%%%%%%%%%%%%%%%%%%%%%%%%%%%%%%%%%%%%%%%%%%%%%%%
% Conclusion + References  -- 1 page.
\appendix

\subsection{Association Probability filter for $\beta_t^m$}
\label{apdx:association_filter_pda}

For the sake of notational simplicity, we will work with
one-dimensional observations, i.e., $h: \Re^d \rightarrow \Re$.  Extension to higher dimensions is straightforward. 

Rewrite~\eqref{eqn:observ_model_seperate} in the vector form:
\begin{equation}
\ud \underline{Z}_t = \underline{\chi}(A_t) h(X_t)\ud t+ \ud \underline{W}_t,\nonumber
%\label{eqn:observ_model_vector}
\end{equation}
where $\underline{Z}_t := (Z_t^1,\hdots,Z_t^M)$, $\underline{\chi}(A_t) := (\chi_t^1,...,\chi_t^M)$,
$\chi_t^m = 1_{[A_t=m]}$, and
$\underline{W}_t:=(W_t^1,...,W_t^M)$. The transition
intensity matrix for the jump Markov process $A_t$ is denoted
as $\Lambda$ with
\begin{equation*}
\Lambda_{mm'} = \begin{dcases*}
        -Mq  &  if $m=m'$\\
        q &  if $m\neq m'$.
        \end{dcases*}
\end{equation*}
Denote $\mathcal{X}_t:= \sigma(X_s:s \leq t)$, and
$\mathcal{C}_t := \mathcal{X}_t \vee \UZ_t$. The
derivation is based on the property of the conditional
expectation:
\begin{equation}
\E[\varphi(A_t)|\UZ_t] = \E[\E[\varphi(A_t)|\mathcal{C}_t]|\UZ_t],\nonumber
\end{equation}
where $\varphi: \{0,\hdots,M\} \rightarrow \Re$ is a test function.

The SDE for evolution of
$\E[\varphi(A_t)|\mathcal{C}_t]$ is described by
the standard Wonham filter:
\begin{align}
\E[\varphi(A_t)|\mathcal{C}_t] = \E[\varphi(A_0)] + \int_0^t \E[\mathcal{L} \varphi(A_s)|\mathcal{C}_s]\ud s+  \int_0^t \E[(\chi(A_s) - \underline{\beta}_s)^T h \varphi(A_s) |\mathcal{C}_s] (\ud \underline{Z}_s  - \underline{\beta}_s h(X_s) \ud s),
\label{eqn:proof_pda_wonham}
\end{align}
where $\underline{\beta}_t := (\beta_t^1,\hdots,\beta_t^M)$ and $\mathcal{L}$ is the forward operator for the Markov process $A_t$ which is defined as:
\begin{equation*}
\mathcal{L} \varphi(m) = \sum_{m'=0}^M \Lambda_{mm'} \varphi(m').
\end{equation*}
Specifically, by choosing 
\[
\varphi(j) = \begin{dcases*}
        1  &  if $j=m$,\\
        0 &  otherwise,
        \end{dcases*}
\]
and taking $\E[\cdot|\UZ_t]$
of~\eqref{eqn:proof_pda_wonham} we obtain:
\begin{equation}
\beta_t^m = \beta_0^m + \int_0^t q\left[1 - (M+1)\beta_s^m\right]\ud s + \int_0^t (\chi(m) - \underline{\beta}_s)^T \hat{h} \beta_s^m (\ud \underline{Z}_s  - \underline{\beta}_s \hat{h} \ud s), \quad m=0,\hdots,M, \nonumber
\end{equation}
which is simplified to obtain~\eqref{eqn:filter_for_beta_nonlinear}-\eqref{eqn:filter_beta_0}.

\subsection{An Alternate Derivation of~\eqref{eqn:filter_for_beta_nonlinear}-\eqref{eqn:filter_beta_0}}
\label{apdx:discretized_assoc_filter}

The aim of this section is to formally derive the update part
of the continuous-time
filter~\eqref{eqn:filter_for_beta_nonlinear}-\eqref{eqn:filter_beta_0} by taking a
continuous-time limit of the discrete-time algorithm for
evaluation of association probability. The procedure for taking
the limit is identical to the original proof of the K-S equation
in~\cite{Kushner64SIAM}, and the proof in Sec~$6.8$ of the classical
text~\cite{jazwinski70}.   

At time $t$, we have $M$ observations $\Delta \underline{Z}_t =
(\Delta Z_t^1,\Delta Z_t^2,...,\Delta Z_t^M)$, at most one of which
originates from the target. The discrete-time filter for
association probability is obtained by using the Bayes' rule:
\begin{equation}
\Prob\{[A_t = m]|\UZ_t,\Delta \underline{Z}_t\} = \frac{\Prob\{\Delta \underline{Z}_t|[A_t = m]\}\Prob\{[A_t = m]|\UZ_t\}}{\sum_{j=0}^M\Prob\{\Delta \underline{Z}_t|[A_t = j]\}\Prob\{[A_t = j]|\UZ_t\}}.\label{eqn:discret_bayes}
\end{equation}
This equation appears also in the derivation of the classical
discrete-time PDA algorithm  (equation~(6-32)
in~\cite{Bar-Shalom_book_88}).  The conditional probabilities are
evaluated as: for $m =1,\hdots,M$,
\begin{equation}
\Prob\{\Delta \underline{Z}_t|[A_t =m]\} =\Prob\{\Delta Z_t^m|[A_t=m]\} \prod_{j\neq m} \Prob_0\{\Delta Z_t^j\}, \label{eqn:Z_t_m}
\end{equation}
and for $m = 0$,
\begin{equation}
\Prob\{\Delta \underline{Z}_t|[A_t = 0]\} = \prod_{j=1}^M \Prob_0\{\Delta Z_t^j\},
\label{eqn:Z_t_0}  
\end{equation}             
where $\Prob\{\Delta Z_t^m|[A_t=m]\} = \frac{1}{(2\pi)^{\frac{s}{2}}\sqrt{\Delta t}}
\int_{\Re^d} \exp\left[-\frac{|\Delta Z_t^m - h(x)\Delta t|^2}{2\Delta
    t}\right]p(x,t)\ud x$ and $\Prob_0\{\Delta Z_t^j\} :=
 \frac{1}{(2\pi)^{\frac{s}{2}}\sqrt{\Delta t}} \exp\left[-\frac{|\Delta
    Z_t^j|^2}{2\Delta t}\right]$.  The formula for $\Prob_0$ is based on the
assumption of the Gaussian clutter model (see Sec.~\ref{sec:prob}).  This is different from the more standard 
model of clutter, popular in discrete-time settings, whereby the
independent observations are uniformly and independently distributed
in a ``coverage area'' $V$ (see~\cite{Bar-Shalom_IEEE_CSM,Bar-Shalom_book_88}). 

Denote $\beta_t^m = \Prob\{[A_t=m]|\UZ_t\}$, the increment in
the observation update step (see Sec~$6.8$
in~\cite{jazwinski70}) is given by
\begin{equation}
\Delta \beta_t^m:= \Prob\{[A_t=m]|\UZ_t,\Delta Z_t\} - \Prob\{[A_t=m]|\UZ_t\}.\label{eqn:def_dbeta}
\end{equation}
Using~\eqref{eqn:discret_bayes} and~\eqref{eqn:def_dbeta}, we
have:
\begin{equation}
\Delta \beta_t^m = E^m(\Delta t,\Delta \underline{Z}_t)\beta_t^m  - \beta_t^m ,\label{eqn:cal_dbeta}
\end{equation}
where 
\begin{equation}
E^m(\Delta t,\Delta \underline{Z}_t) := \frac{\Prob\{[A_t = m]|\UZ_t,\Delta \underline{Z}_t\}}{\Prob\{[A_t = m]|\UZ_t\}} = \frac{\Prob\{\Delta \underline{Z}_t|[A_t = m]}{\sum_{j=0}^M\Prob\{\Delta \underline{Z}_t|[A_t = j]\}\Prob\{[A_t = j]|\UZ_t\}}.\nonumber
\end{equation}
Using It$\hat{\text{o}}$'s rule,
\begin{equation}
(\Delta Z_t^j)^T \Delta Z_t^k = \begin{dcases*}
        \Delta t,&\; {\text{if}} j = k, \\
        0, &\;\text{otherwise,}
        \end{dcases*}\nonumber
\end{equation}
we have for $m = 1,\hdots,M$,
\begin{equation}
E^m(\Delta t,\Delta \underline{Z}_t)= \frac{\int \exp\left[h^T(x) \Delta Z_t^m - \frac{1}{2}h^T(x)h(x)\Delta t\right]p(x,t)\ud x}{\sum_{l=1}^M \beta_t^l \int \exp\left[(h^T(x) \Delta Z_t^l - \frac{1}{2}h^T(x)h(x)\Delta t\right]p(x,t)\ud x + \beta_t^0},\nonumber
\end{equation}
and for $m=0$,
\begin{equation}
E^0(\Delta t,\Delta \underline{Z}_t) = \frac{1}{\sum_{l=1}^M \beta_t^l \int \exp\left[(h^T(x) \Delta Z_t^l - \frac{1}{2}h^T(x)h(x)\Delta t\right]p(x,t)\ud x + \beta_t^0}.\nonumber
\end{equation}
Following the derivation procedure in~\cite{jazwinski70}, we expand $E^m(\Delta t, \Delta \underline{Z}_t)$ as a multivariate
series about $(0,\underline{0})$:
\begin{align}
E^m(\Delta t, \Delta \underline{Z}_t) = E^m(0,\underline{0}) + E^m_{\Delta t}(0,\underline{0})\Delta t + \sum_{j=1}^M  E^m_{\Delta Z_t^j} (0,\underline{0}) \Delta Z_t^j + \frac{1}{2}\sum_{j,k=1}^M (\Delta Z_t^j)^T E^m_{\Delta Z_t^j,\Delta Z_t^k}(0,\underline{0}) \Delta Z_t^k + o(\Delta t).\nonumber
\end{align}
By direct evaluation, we obtain: for $m = 1,\hdots,M$,
\begin{align}
E^m(0,\underline{0}) = 1,& \quad E^m_{\Delta t}(0,\underline{0}) = -\frac{1}{2} \beta_t^0\widehat{|h^2|},\nonumber\\
E^m_{\Delta Z_t^j} (0,\underline{0}) = -\beta_t^j \hat{h}^T ,& \quad E^m_{\Delta Z_t^j,\Delta Z_t^j}(0,\underline{0}) = 2(\beta_t^j)^2 |\hat{h}|^2 - \beta_t^j \widehat{|h^2|},\quad  j \neq m, \nonumber\\
E^m_{\Delta Z_t^m} (0,\underline{0}) = (1-\beta_t^m) \hat{h}^T,& \quad
E^m_{\Delta Z_t^m,\Delta Z_t^m}(0,\underline{0}) = (1-\beta_t^m)\widehat{|h^2|} -  2\beta_t^m(1-\beta_t^m)|\hat{h}|^2 ,\nonumber
\end{align}
and for $m = 0$,
\begin{align}
E^0(0,\underline{0}) = 1,&\quad E^0_{\Delta t}(0,\underline{0}) = \frac{1}{2}(1-\beta_t^0) \widehat{|h^2|},\nonumber\\
E^0_{\Delta Z_t^j} (0,\underline{0}) = -\beta_t^j \hat{h}^T,&\quad E^0_{\Delta Z_t^j,\Delta Z_t^j}(0,\underline{0}) = 2(\beta_t^j)^2 |\hat{h}|^2 - \beta_t^j \widehat{|h^2|},\quad  j = 1,\hdots,M, \nonumber
\end{align}
where $\hat{h} := \E[h(X_t)|\UZ_t]$, $|\hat{h}|^2 := \hat{h}^T \hat{h}$ and
$\widehat{|h^2|}:=\E[|h(X_t)|^2|\UZ_t]$. 

By using It$\hat{\text{o}}$'s rule again, we obtain the simplified expression: for $m = 1,\hdots,M$,
\begin{equation}
E^m(\Delta t,\Delta \underline{Z}_t) =  1 + \hat{h}^T \left(\Delta Z_t^m - \sum_{j=1}^M \beta_t^j\Delta Z_t^j\right) + |\hat{h}|^2 \left(\sum_{j=1}^M (\beta_t^j)^2 - \beta_t^m\right)\Delta t,\label{eqn:Em_approx_m}
\end{equation}
and for $m = 0$,
\begin{equation}
E^0(\Delta t,\Delta \underline{Z}_t) = 1 - \hat{h}^T \sum_{m=1}^M \beta_t^m \Delta Z_t^m + |\hat{h}|^2 \sum_{m=1}^M (\beta_t^m)^2 \Delta t.
\label{eqn:Em_approx_0}
\end{equation}

Substituting~\eqref{eqn:Em_approx_m}-\eqref{eqn:Em_approx_0} into~\eqref{eqn:cal_dbeta}, we
obtain the expression for $\Delta \beta_t^m$. This is identical to the observation update part of the continuous-time
filter~\eqref{eqn:filter_for_beta_nonlinear}-\eqref{eqn:filter_beta_0}.
\begin{remark}
\label{rem:remark-DT-DAP}
In a discrete-time implementation, one can use~\eqref{eqn:discret_bayes}-\eqref{eqn:Z_t_0} to compute
the association probability. In the evaluation of~\eqref{eqn:Z_t_m}, $\Prob\{\Delta Z_t^m|[A_t=m]\}$ is approximated by using particles:
\begin{equation}
\Prob\{\Delta Z_t^m|[A_t=m]\}\approx \frac{1}{N}\frac{1}{(2\pi)^{\frac{s}{2}}\sqrt{\Delta t}} \sum_{i=1}^N \exp\left[-\frac{|\Delta Z_t^m - h(X_t^i)\Delta t|^2}{2\Delta t}\right].\nonumber
\end{equation}
In this case, the formula for Gaussian clutter observation density, $\Prob_0\{\Delta
Z_t^j\}$ in~\eqref{eqn:Z_t_m}-\eqref{eqn:Z_t_0}, can be replaced another commonly used
clutter observation density model, e.g., uniform distribution in a ``coverage
area'' $V$~\cite{Bar-Shalom_IEEE_CSM,Bar-Shalom_book_88}.  One can also
incorporate gating and detection probability
into~\eqref{eqn:discret_bayes}-\eqref{eqn:Z_t_0} as in the classical
literature. 

The important point is that once the association probability is computed, the filtering
equation for the state process -- the PDA-FPF algorithm -- 
does {\em not} depend upon the details of the clutter model.  The
Gaussian clutter model assumption is primarily used here in the
derivation of the continuous-time
filter~\eqref{eqn:filter_for_beta_nonlinear}-\eqref{eqn:filter_beta_0} for computing the
association probability.        
\qed
\end{remark}

\subsection{Consistency proof of $p$ and $p^\ast$}
\label{apdx:consistency}

\noindent {\bf Evolution equation for $p^\ast$:} Recall
$\mathcal{\underline{Z}}_t := \sigma(\underline{Z}_s:s\le t)$,
$\mathcal{A}_t:=\sigma(A_s:s\le t)$. We denote $\mathcal{C}_t
:= \mathcal{A}_t \vee \mathcal{\underline{Z}}_t$.  The
derivation is based on the property of the conditional
expectation:
\begin{equation*}
{\sf E} [ \varphi(X_t)|\mathcal{\underline{Z}}_t] = {\sf E} \left[
  {\sf E} [ \varphi(X_t)|\mathcal{{C}}_t] |
  \mathcal{\underline{Z}}_t \right].
\end{equation*}
The sde for evolution of  ${\sf E} [
\varphi(X_t)|\mathcal{{C}}_t]$ is described by the standard
nonlinear filter with innovation error, $\sum_{m=1}^M \chi_t^m (\ud
Z_t^m - \hat{h} \ud t)$, where $\chi_t^m = 1_{[A_t=m]}$:
\begin{align*}
{\sf E} [\varphi(X_t)|\mathcal{C}_t] &= {\sf E} [
\varphi(X_0)] +  \int_{0}^t {\sf E} [
\clL \varphi(X_s)|\mathcal{C}_s] \ud s + \sum_{m=1}^M  \int_{0}^t {\sf E} [(
h-\hat{h})^T \varphi(X_s) |\mathcal{C}_s]  \chi_s^m (\ud
Z_s^m - \hat{h} \ud s),
\end{align*}
where $\clL$ denotes the Kolmogorov's backward operator for the diffusion~\eqref {eqn:Signal_Process} (the adjoint of  $\clL^\dagger$).

Taking $  {\sf E} [\cdot|\mathcal{\underline{Z}}_t]$ gives the desired
result because ${\sf E} [\chi_s^m |\mathcal{\underline{Z}}_s] = {\sf
  P}\{[A_s = m]|\mathcal{\underline{Z}}_s\} = \beta_s^m$.

\smallskip
\noindent {\bf Evolution equation for $p$:}
We express the feedback particle filter~\eqref{eqn:PDA_FPF} as:
\[
\ud X_t^i = a(X^i_t)\ud t + \sigma_B \ud B_t^i +  \v(X_t^i,t) \sum_{m=1}^{M}
\beta_t^m \ud Z^m_t + u(X_t^i,t) \ud t,
\]
where
\begin{align}
u(x,t):= - \sum_{m=1}^{M} \beta_t^m \v(x,t) \left[ \frac{\beta_t^m }{2} h +
  (1 - \frac{\beta_t^m }{2}) \hat{h} \right]   +
 \sum_{m=1}^{M} (\beta_t^m)^2 \Omega(x,t), \label{eqn:PDA_u_def}
\end{align}
and $\Omega(x,t) = (\Omega_1,\Omega_2,\hdots,\Omega_d)$ is the Wong-Zakai correction term for~\eqref{eqn:PDA_FPF}:
\begin{equation*}
\Omega_l(x,t) := \frac{1}{2}\sum_{k=1}^d \sum_{j=1}^s \v_{kj}(x,t) \frac{\partial \v_{lj}}{\partial x_k}(x,t),\;\; \text{for} \;\; l \in \{1,\hdots,d\}.
\end{equation*}

The evolution equation for $p$ now follows:
\begin{align}
\ud p = \clL^\dagger p \ud t  - \nabla \cdot (pu) \ud t - \sum_{m=1}^{\M} \beta_t^m \nabla \cdot (p\v) \ud Z_t^m + \sum_{m=1}^{\M} (\beta_t^m)^2 \left(\frac{1}{2}\sum_{l,k=1}^d \frac{\partial^2}{\partial x_l \partial x_k}(p[\v \v^T]_{lk})\right)\ud t.
%\frac{\partial}{\partial x}( u p ) \ud t
% &-\sum_{j=1}^{m_t} \left(\beta_t^2(j) \frac{\partial}{\partial x}( u p ) + (\beta_t(j)+\beta_t^2(j)) \frac{\partial}{\partial x}(p\v\hat{h})\right) \ud t\nonumber\\
%+  \frac{\sigma_W^2}{2} \sum_{j=1}^{M}
%(\beta_t^m)^2\frac{\partial^2}{\partial x^2}\left( p \v^2 \right) \ud t \nonumber\\
%& - \frac{\partial}{\partial x}\left( \v
 % p \right) \sum_{m=1}^{M} \beta_t^m \ud Z_t^m.
\label{eqn:mod_PDA_FPK}
\end{align}

\smallskip

\noindent {\bf Proof of consistency:} The proof follows closely
the consistency proof for the multivariable feedback particle filter (see
Appendix.~A in~\cite{taoyang_cdc12}).  If $\v$ solves the
E-L BVP then
\begin{equation}
\nabla \cdot (p \v) = -(h-\hat{h})^Tp.\label{eqn:pv}
\end{equation}
On multiplying both sides
of~\eqref{eqn:PDA_u_def} by $-p$ and simplifying (by using~\eqref{eqn:pv}),
we obtain
\begin{align}
-up &= \sum_{m=1}^\M \frac{1}{2} (\beta_t^m)^2 \v (h-\hat{h}) p + \sum_{m=1}^\M \beta_t^m p \v \hat{h} - \sum_{m=1}^M (\beta_t^m)^2 p \Omega(x,t)\nonumber\\
&= -\sum_{m=1}^M \frac{1}{2}(\beta_t^m)^2 \v \left[\nabla \cdot(p\v)\right]^T - \sum_{m=1}^M (\beta_t^m)^2 p \Omega + \sum_{m=1}^M \beta_t^m p \v \hat{h},\label{eqn:multi_up}
\end{align}
where~\eqref{eqn:pv} is used to obtain the second equality. Denoting $E:= \frac{1}{2} \v [\nabla \cdot (p\v)]^T$, a direct calculation shows that
\begin{equation*}
E_{l} +  \Omega_{l}p = \frac{1}{2}\sum_{k=1}^d \frac{\partial }{\partial x_k}
\left( p[\v \v^T]_{lk} \right).
\end{equation*}
Substituting this in~\eqref{eqn:multi_up}, on taking the divergence of both sides, we obtain:
\begin{align}
-\nabla \cdot (pu) +\sum_{m=1}^M (\beta_t^m)^2 \frac{1}{2}\sum_{l,k=1}^d \frac{\partial^2}{\partial x_l \partial x_{k}} &\left( p [\v \v^T]_{lk} \right) = \sum_{m=1}^M \beta_t^m \nabla \cdot (p\v) \hat{h}.
\label{eqn:FPK_23}
\end{align}

Using~\eqref{eqn:pv} and~\eqref{eqn:FPK_23}
in the forward equation~\eqref{eqn:mod_PDA_FPK}, we obtain:
\begin{equation}
\ud p = \clL^\dagger p \ud t + \sum_{j=1}^{M} \beta_t^m (h-\hat{h})^T (\ud Z_t^m- \hat{h}\ud t)p.\nonumber
%\label{eqn:PDA_p_FPK}
\end{equation}
This is precisely the SDE~\eqref{eqn:PDA_KS}, as desired.

\subsection{Joint probabilistic data association-feedback particle filter}
\label{apdx:jpda_algo}

%Since the association events are independent, we have:
%\[
%\pi_t^{\ugamma} = \prod_{m=1}^M \P\{[\alpha_t^m = \gamma^m]|\UZ_t\}.
%\]

By repeating the steps of the proof in Appendix.~\ref{apdx:association_filter_pda}, the filter for the joint association probability $\pi_t^{\ugamma}$ (defined as in~\eqref{eqn:def_joint_prob}) is of the following form:
\begin{align}
\ud \pi_t^{\ugamma} &= q \left[1 -M!\, \pi_t^{\ugamma}\right] \ud t + \pi_t^{\ugamma} \sum_{m=1}^M \left[\widehat{h(X_t^{\gamma^m})} - \widehat{\tilde{H}_t^m}\right]^T \ud Z_t^m - \pi_t^{\ugamma} \sum_{m=1}^M \left[\widehat{h(X_t^{\gamma^m})}^T \widehat{\tilde{H}_t^m} - \left|\widehat{\tilde{H}_t^m}\right|^2\right] \ud t,\label{eqn:wonham_jpdaf}
\end{align}
where $\tilde{H}_t = (\tilde{H}_t^1,\hdots,\tilde{H}_t^M) := \sum_{\{\ugamma \in \Pscr(\Nscr)\}} \pi_t^{\ugamma} H(\ugamma)$, $H(\ugamma) := \left(h(X_t^{\gamma^1}),\hdots,h(X_t^{\gamma^M})\right)$, $\widehat{h(X_t^{\gamma^m})}:=\E[h(X_t^{\gamma^m})|\UZ_t]$, $\widehat{\tilde{H}_t^m}:=\E[\tilde{H}_t^m|\UZ_t]$ and $\left|\widehat{\tilde{H}_t^m}\right|^2 := \widehat{\tilde{H}_t^m}^T \widehat{\tilde{H}_t^m}$. All the conditional expectation are approximated by using particles. %The derivation of~\eqref{eqn:wonham_jpdaf} is a direct extension of the proof in Appendix~\ref{apdx:association_filter_pda} and omitted here on account of space

Denote $\beta_t^{m,n}$ to be the conditional probability that the $\mth$ observation originates from the $\nth$ target:
\begin{equation}
\beta_t^{m,n} := \P\{[\alpha_t^m = n]|\UZ_t\}.\nonumber
%\label{eqn:marginal_assoc}
\end{equation}
It is obtained by using joint association probability $\pi_t^{\ugamma}$:
\begin{equation}
\beta_t^{m,n} = \sum_{\{\ugamma \in \Pscr(\Mscr): \gamma^m = n\}} \pi_t^{\ugamma}. \nonumber
%\label{eqn:ProbAsso}
\end{equation}

The JPDA-FPF for the $\nth$ target is given by the following controlled system:
\begin{align}
\ud X_t^{i;n} = a^n(X^{i;n}_t)\ud t + \ud B_t^{i;n} + \sum_{m=1}^{\M} \beta_t^{m,n} \, \v^n(X_t^{i;n},t) \circ \ud \Inov^{i,m;n}_t,\label{eqn:JPDA_FPF}
\end{align}
where $X^{i;n}_t \in \Re^d$ is the state for the $i^{\text{th}}$
particle at time $t$ for the $\nth$ target, $\{B^{i;n}_t\}$ are mutually independent standard Wiener processes, $\Inov^{i,m;n}_t$ is a modified form of the
$\emph{innovation process}$,
\begin{equation}
\ud \Inov^{i,m;n}_t := \ud Z_t^m - \left[\frac{\beta_t^{m,n}}{2} h(X_t^{i;n}) + \left(1-\frac{\beta_t^{m,n}}{2}\right)\hat{h}^n\right]\ud t,\label{eqn:jpda_inov}
\end{equation}
where $\hat{h}^n := {\sf E} [h(X_t^{i;n})|\UZ_t] = \int_{\Re^d} h(x)
p_n(x,t) \ud x \approx \frac{1}{N} \sum_{j=1}^N h(X_t^{j;n})$ and $p_n(x,t)$ denotes the conditional distribution of
    $X^{i;n}_t$ given $\UZ_t$.

The gain function $\v^n=[\nabla \phi_1^n,\hdots,\nabla \phi_s^n]$ is a solution of a certain E-L BVP: for $j\in\{1,\hdots,s\}$,
\begin{equation}
\nabla \cdot (p_n \nabla \phi_j^n) = - (h_j-\hat{h}_j^n) p_n.\label{eqn:jpda_bvp}
\end{equation}
%In a numerical implementation, we approximate $\nabla \phi_j^n$ using particles:
%\begin{equation}
%\nabla \phi_j^n \approx \frac{1}{N} \sum_{i=1}^N X_t^{i;n}\left( h_j(X_t^i) - \hat{h}_j^n\right),\nonumber
%\end{equation}
%where $\hat{h}_j^n \approx \frac{1}{N}h_j(X_t^{i;n})$.

\begin{remark}
The corresponding ODE model for the JPDA-FPF~\eqref{eqn:JPDA_FPF}-\eqref{eqn:jpda_inov} is given by:
\begin{align}
&\frac{\ud X_t^{i;n}}{\ud t} = a^n(X^{i;n}_t) + \dot{B}_t^{i;n} +\sum_{m=1}^{\M} \beta_t^{m,n} \, \v^n(X_t^{i;n},t) \left( Y_t^m - \left[\frac{\beta_t^{m,n}}{2}h(X_t^{i;n}) + \left(1 - \frac{\beta_t^{m,n}}{2}\right)\hat{h}^n\right] \right),\nonumber
%\label{eqn:ode_jpda_fpf}
\end{align}
where $\{\dot{B}_t^{i;n}\}$ are independent white noise processes and $Y_t^m \doteq \frac{\ud Z_t^m}{\ud t}$. The gain function $\v^n$ is again a solution of the E-L BVP~\eqref{eqn:jpda_bvp}, and may be approximated by using~\eqref{eqn:prelim_const_gain}.
\qed
\end{remark}

\begin{remark}
In a discrete-time implementation, one can use the following heuristic to obtain association probability:
\begin{align}
    \pi_t^{\ugamma} &= \prod_{m=1}^M \P\left\{[\alpha_t^m = \gamma^m]|\UZ_t\right\}\; (\text{Assoc. events are independent})\nonumber\\
    &\propto \prod_{m = 1}^\N \Prob\left\{
    Y_t^{m}\,\Big|\, [\alpha_t^{m} = \gamma^m] \right\}\nonumber\\
    &= \frac{1}{(2\pi)^{\frac{s}{2}}} \prod_{m = 1}^\N \left(\int_{\R^d} \exp \left[
    - \frac{1}{2}| Y_t^{m}-h(x) |^2 \right] p_{\gamma^m}(x, t) \ud x \right) \nonumber\\
    &\approx \frac{1}{(2\pi)^{\frac{s}{2}}} \prod_{m = 1}^\N \left(
    \frac{1}{N}\sum_{i=1}^N \exp \left[ - \frac{1}{2}\left| Y_t^{m}-h(\particleX_t^{i;\gamma^m}) \right|^2 \right]\right).\nonumber
    %\label{eqn:approx_joint_assoc}
\end{align}
%
%The discrete-time version algorithm for the association probability filter is useful for applications where observations are made at discrete sampling times.
\qed
\end{remark}

%\subsection{Association Probability filter for $\pi_t^m$}
%\label{apdx:association_filter_jpda}
%
%The derivation follows closely the derivation in
%Appendix~\ref{apdx:association_filter_pda}. Note that the
%observation model is described
%by~\eqref{eqn:Two_Target_Observ_Process}. Denote
%$\underline{h}(\underline{X}_t) := (h(X_t^1),h(X_t^2))^T$,
%$\widehat{\Psi}_t := \sum_{m=1}^2 \pi_t^m \Psi(m)$ and
%$\underline{\phi}_t = (\phi_t^1,\phi_t^2)^T:=\widehat{\Psi}_t
%\underline{h}(\underline{X}_t)$. The Wonham filter is given by:
%\begin{align}
%&\E[\underline{\varphi}(A_t)|\mathcal{C}_t] = \E[\underline{\varphi}(A_0)] + \int_0^t \E[\Lambda \underline{\varphi}(A_s)|\mathcal{C}_s]\ud s\nonumber\\
%&+ \sum_{m=1}^2 \int_0^t \E[(D_s^m - \phi_s^m I)^T \underline{\varphi}(A_s)|\mathcal{C}_s] (\ud Z_s^m - \phi_s^m \ud s),
%\label{eqn:proof_jpda_wonham}
%\end{align}
%where $D_t^m$ is a $2\times 2$ diagonal matrix where
%$(D_t^m)_{ii}$ is the $m^{\text{th}}$ entry of the vector
%$\Psi(i) \underline{h}(\underline{X}_t)$.
%
%
%Taking $\E[\cdot|\mathcal{Z}_t]$
%of~\eqref{eqn:proof_jpda_wonham} gives the desired result.
%

%%%%%%%%%%%%%%%%%%%%%%%%%%%%%%%%%%%%%%%%%%%%%%%%%%%%%%%%%%%%%%%%%%%%%%%%%%%%%%
\bibliographystyle{IEEEtran}
%\bibliography{ACCPF,ACCJPDAFPF,refmtt,fpfbib,JPDAF_journal,JPDAFPF_Tran}
\bibliography{fpfbib,JPDAFPF_Tran}

%\bibliography{ref}
%%%%%%%%%%%%%%%%%%%%%%%%%%%%%%%%%%%%%%%%%%%%%%%%%%%%%%%%%%%%%%%%%%%%%%%%%%%%%%

\end{document}